\documentclass[MSNbibl,number,citesort,seceqn,dvips]{arxbj}

\aid{0}
\volume{19}
\issue{5A}
\pubyear{2013}
\firstpage{2098}
\lastpage{2119}
\doi{10.3150/12-BEJ444} 

\makeatletter

\newtheorem{theorem}{Theorem}[section]
\newremark{example}{Example}[section]
\newtheorem{lemma}{Lemma}[section]
\newremark{remark}{Remark}[section]
\newtheorem{corollary}{Corollary}[section]
\makeatother

\begin{document}
\begin{frontmatter}

\title{Marked empirical processes for non-stationary time series}
\runtitle{Empirical processes}

\begin{aug}
\author[1]{\fnms{Ngai Hang} \snm{Chan}\thanksref{1}\ead[label=e1]{nhchan@sta.cuhk.edu.hk}} \and
\author[2]{\fnms{Rongmao} \snm{Zhang}\corref{}\thanksref{2}\ead[label=e2]{rmzhang@zju.edu.cn}}
\runauthor{N.H. Chan and R. Zhang} 
\address[1]{Department of Statistics, The Chinese University of Hong Kong,
Shatin, NT, Hong Kong.\\ \printead{e1}}
\address[2]{Department of Mathematics, Zhejiang University, Hangzhou, China.
\printead{e2}}
\end{aug}

\received{\smonth{9} \syear{2011}}
\revised{\smonth{3} \syear{2012}}

%
\begin{abstract}
Consider a first-order autoregressive process $X_i=\beta
X_{i-1}+\varepsilon_i,$ where $\varepsilon_i=G(\eta_i, \eta_{i-1},
\ldots)$
and $\eta_i, i\in\mathbb{Z} $ are i.i.d. random variables. Motivated
by two important issues for the inference
of this model, namely, the quantile inference for $H_0\dvt\beta=1$,
and the goodness-of-fit for the unit root model, the notion of the
marked empirical process
$\alpha_n(x)=\frac{1}{ n}\sum_{i=1}^{n} g(X_i/a_n) I(\varepsilon_i\leq
x), x\in\mathbb{R}$
is investigated in this paper. Herein, $g(\cdot)$ is a continuous
function on $\mathbb{R}$ and $\{a_n\}$ is a sequence of
self-normalizing constants.
As the innovation $\{\varepsilon_i\}$ is usually not observable, the
residual marked empirical process
$\hat\alpha_n(x)=\frac{1}{ n}\sum_{i=1}^{n} g(X_i/a_n) I(\hat
\varepsilon_i\leq x), \allowbreak x\in\mathbb{R},$
is considered instead, where $\hat\varepsilon_i=X_i-\hat\beta
X_{i-1}$ and $ \hat\beta$ is a consistent estimate of $\beta.$ In
particular, via the martingale decomposition of stationary process and
the stochastic integral result of Jakubowski (\textit{Ann. Probab.} \textbf{24}
(1996) 2141--2153),
the limit distributions of $\alpha_n(x)$ and $\hat\alpha_n(x)$ are
established when $\{\varepsilon_i\}$ is a short-memory process.
Furthermore, by virtue of the results of Wu (\textit{Bernoulli} \textbf{95}
(2003) 809--831) and Ho and Hsing
(\textit{Ann. Statist.} \textbf{24} (1996) 992--1024) of empirical process and
the integral result of
Mikosch and Norvai\v{s}a (\textit{Bernoulli} \textbf{6} (2000) 401--434) and
Young (\textit{Acta Math.} \textbf{67} (1936) 251--282), the limit
distributions of $\alpha_n(x)$ and $\hat\alpha_n(x)$ are also
derived when $\{\varepsilon_i\}$ is a long-memory process.

\end{abstract}

%

\end{frontmatter}
%

\section{Introduction}\label{sec1}

Consider the autoregressive (AR) model
%
\begin{equation}
\label{1.1}X_i=\beta X_{i-1}+\varepsilon_i,
\end{equation}
where $X_0$ is given and $ \varepsilon_i=G(\eta_i, \eta_{i-1},
\ldots, )$ is such that $\mathrm{E}\varepsilon_i=0$ (when it exists)
and $\{\eta_i\}$ is a sequence of i.i.d. random variables. It is known
that when the tail of $\varepsilon_i$ is heavy, the quantile estimate
of $\beta$ performs better than the least squares estimate (LSE). Even
under the Gaussian setting, Zou and Yuan \cite{zy08} proved that a composite
quantile estimate can be as efficient as
the maximum likelihood estimate (MLE). As a result, the quantile
estimate provides a good alternate to the LSE. The first issue pursued
in this
paper is to study the asymptotic properties of the quantile estimate
for model (\ref{1.1}) with both long and short-memory
innovations when $\beta=1$.

A second motivation of this paper is to consider the goodness-of-fit
issue for model (\ref{1.1}). Empirical processes and goodness-of-fit
tests in the i.i.d. case have long been a vibrant research
topic in statistics, see, for example, the succinct monograph of del
Barrio, Deheuvels and van de Geer \cite{bdg07}, the proceeding of Gaenssler
and R\'{e}v\'{e}sz \cite{gr76} and the references therein. Recently, there
have also been developments on model checking using goodness-of-fit
ideas for dependent data. For example, Bai \cite{b03} applied a
Rosenblatt-transform to test the conditional distribution of
$\varepsilon_1$ under condition on $\{\varepsilon_i, i\leq0\}$,
Escanciano \cite{e06} and Hong and Lee \cite{hl03} used a generalized spectral
method to check the model fitness, Koul and Ling \cite{kl06}
considered the
Kolmogorov--Smirnov (\mbox{K--S}) statistics of empirical process for GARCH
model and Chan and Ling \cite{cl08} generalized the K--S test to long-memory
time series. It should be noted that Chan and Ling
\cite{cl08} only made use of the marginal distribution information of
$\varepsilon_1$. For long-memory dependence, using only the marginal
distribution information may reduce the test power and lead to
incorrect conclusions. An alternative statistic which increases
the power and takes into account of the dependent information is
therefore required. Recently, some progresses have been made on this
issue. For example, Woodridge \cite{w90} and Escanciano \cite{e07}
proposed the
statistic
$\sum_{i=1}^{n}g(X_{i-1})\varepsilon_i$ for a measurable weighted
function $g(\cdot)$. Stute, Xu and Zhu \cite{sxz08} used $\sum_{i=1}^{n}g(X_{i-1})(I(\varepsilon_i\leq x)-F(x))$
to test the validity of a model for independent data. The idea of
Stute, Xu and Zhu was used by Escanciano \cite{e10} to check the joint
specification of conditional mean and variance of a GARCH-type model.

It turns out that the key idea of studying these two issues lies in
analyzing the asymptotic property of
%
\begin{equation}
\label{1.2} \sum_{i=1}^{n}g(X_{i-1})
\bigl[I(\varepsilon_i<x)-F(x) \bigr], \qquad x\in\mathbb{R},
\end{equation}
where $g(\cdot)$ is a measurable weighted continuous function on
$\mathbb{R}.$
Note that if $X_n$ is a unit root model, then under some regularity conditions,
there exists a constant sequence $\{a_n\}$ such that
$X_{[nt]}/a_n\stackrel{\mathrm{f.d.d.}}{\longrightarrow}
\xi(t)$ for some random process $\xi(t)$, where $\stackrel
{\mathrm{f.d.d.}}{\longrightarrow}$ denotes the weak convergence of
finite-dimensional distributions. This leads us to replace the
statistic of Stute, Xu and Zhu \cite{sxz08} by
%
\begin{equation}
\label{1.3} \alpha_n(x)=\sum_{i=1}^{n}g(X_{i-1}/a_n)
\bigl(I(\varepsilon_i\leq x)-F(x) \bigr).
\end{equation}
%
Observe that $\alpha_n(x)$ is a general form of (\ref{1.2}) and its
limit behavior offers important insight
in studying the two aforementioned issues.

Specifically, let $\{\eta_i'\}$ be an i.i.d. copy of $\{\eta_i\}$,
$\mathcal{F}_i$ be the $\sigma$-field generated by $\{\eta_t,
t\leq i\}$, that is, $\mathcal{F}_i=\sigma(\eta_{i}, \eta_{i-1},
\ldots)$ and
$\mathcal{F}_i^{*}=\sigma(\eta'_i, \eta_{i-1}, \ldots, \eta_{1},
\eta_0, \eta_{-1}, \eta_{-2}\ldots,)$. Let $L^p$ be the space of
random variables $Z$ with $\Vert Z\Vert_p=(\mathrm
{E}|Z|^p)^{1/p}<\infty$.
For simplicity, we also write $\Vert\cdot\Vert_2$ as $\Vert\cdot
\Vert.$ For $j\in
\mathbb{Z}$, define the projection operator
\[
\mathcal{P}_{j}\cdot=\mathrm{E}(\cdot|\mathcal{F}_j)-
\mathrm{E}(\cdot|\mathcal{F}_{j-1})
\]
and define the predictive dependence measure $\theta_p(i)=\Vert
\mathcal
{P}_{0}\varepsilon_i\Vert_p$ as in Wu \cite{w07}. We say that a process
$\{\varepsilon_i\}$ is a short-memory process if $\sum_{i=0}^{\infty
}\theta_p(i)<\infty,$ otherwise, we say it is a long-memory process.

The main purpose of this paper is to consider a unified approach for
the limit of (\ref{1.3}) and the statistic
%
\begin{equation}
\label{1.4} \hat\alpha_n(x)=\sum_{i=1}^{n}g(X_{i-1}/a_n)
\bigl(I(\hat\varepsilon_i\leq x)-F(x) \bigr)
\end{equation}
for model (\ref{1.1}) under non-stationarity
with long and/or short-memory innovations $\{\varepsilon_i\},$ where
$\hat\varepsilon_i=X_i-\hat{\beta}X_{i-1}$ and $ \hat{\beta}$ is
an estimate of $\beta.$
Although Escanciano \cite{e10} has also considered the model checking
problem for dependent data (GARCH), his underlying model still
possesses a martingale structure.
When $\{X_i\}$ is a stationary process adapted to the fields $\mathcal
{F}_i=\sigma(\varepsilon_i, \varepsilon_{i-1}, \ldots)=
\sigma(\eta_i, \eta_{i-1}, \ldots)$ and $\{\varepsilon_i\}$ is a
martingale difference sequence, the central limit theorem (CLT) for
martingale differences can be applied to derive the limit distribution
of (\ref{1.3}). When $X_i$ is a random walk process, it is not clear
how to derive the limit distribution of (\ref{1.3}), especially, if $\{
\varepsilon_i\}$ is a long-memory process. This is because the CLT of
martingale differences cannot be directly used and when $\{\varepsilon_i\}$ is a long-memory process, neither
$g(X_{[nt]}/a_n)$ nor $\sum_{i=1}^{[nt]}(I(\hat\varepsilon_i<x)-F(x))$ can be approximated by a
martingale. In this paper, we first use the result of Jakubowski \cite{j96}
(see also
Protter and Kurtz \cite{kp91}) on the weak convergence of
the stochastic integral to deduce the limit distributions of (\ref
{1.3}) and (\ref{1.4}), when $\{\varepsilon_i\}$ is a short-memory process.
We then combine the results of Wu \cite{w03} and Ho and Hsing \cite
{hh96} on the
expression of empirical process and the integral result of Mikosch and
Norvai\v{s}a \cite{mn00} and Young \cite{y36} to deduce the limit
distributions
of $\alpha_n(x)$ and $\hat\alpha_n(x)$ when $\{\varepsilon_i\}$ is
a long-memory process.

The paper is organized as follows. In Section \ref{sec2}, we consider
the marked
empirical process when $\{\varepsilon_i\}$ is short-memory. Section
\ref{sec3}
considers the case with long-memory innovations. Proofs are given in
Section \ref{sec4}.

\section{Short-memory error processes}\label{sec2}

In this section, we consider the limit distribution for (\ref{1.3})
and (\ref{1.4}) when $\{\varepsilon_i=
G(\eta_i,$ $ \eta_{i-1}, \ldots\}$ is a short-memory process with
mean zero. Let $S_{[nt]}=:S_{n}(t)=\sum_{i=1}^{[nt]}\varepsilon_i$.
According to Theorem~2* of Chapter 7 (see pages 162 and 175) of
Gnedenko and Kolmogorov \cite{gk54},
if $\{\varepsilon_i\}$ are i.i.d. and there exists a sequence $\{a_n\}
$ such that
$S_n(1)/a_n\stackrel{\mathcal{L}}{\longrightarrow} S(1),$
then $S(1)$ is a stable variable, where $\stackrel{\mathcal
{L}}{\longrightarrow}$ denotes the convergence in distribution.
Further, when $\varepsilon_i$ has an infinite variance, it must
satisfy for any $y>0,$
%
\begin{eqnarray}
\label{2.0a}\lim_{x\rightarrow\infty}P\bigl(|\varepsilon_1|\geq xy\bigr)/P\bigl(|
\varepsilon_1|\geq x\bigr)=y^{-\alpha}
\end{eqnarray}
and the normalization
constants $\{a_n\}$ are given by
%
\begin{eqnarray}
\label{2.0b}a_n=\inf \bigl\{x\dvt P\bigl(|\varepsilon_1|\geq
x\bigr)\leq1/n \bigr\} .
\end{eqnarray}
Similar behaviors exist for short-memory processes under certain
regularity conditions, see, for example, Davis and Resnick \cite{dr86}.
Throughout the paper,
we assume (\ref{2.0a}) and (\ref{2.0b}) hold when
$\{\varepsilon_i\}$ has an infinite variance and there exists a
constant sequence $\{a_n\}$
such that \mbox{$S_n(1)/a_n\stackrel{\mathcal{L}}{\longrightarrow} S(1).$}

Let $F_{i}(x|\mathcal{F}_j)=P(\varepsilon_i\leq x|\mathcal{F}_j),$ $
f_{i}(x|\mathcal{F}_j)=F_{i}^{(1)}(x|\mathcal{F}_j)$ be the
conditional distribution
(resp., density) function of $\varepsilon_i$ at $x$ given
$\mathcal{F}_j$ and $f_i$ be the marginal density of $\varepsilon_i$. Let
$ W_n(t, x)=\sum_{i=1}^{[nt]}[I(\varepsilon_i<x)-F(x)]$ and $W(t, x)$
be a rescaled Brownian bridge for fix $t$ and a Brownian motion with
variance $\mu(x)=E\{\sum_{i=0}^{\infty} F_i(x|\mathcal
{F}_0)-F_i(x|\mathcal{F}_0^*)\}^2$ for fix $x.$ Further, let
$\stackrel{S}{\Longrightarrow}$ and $\stackrel{w}{\Longrightarrow}$
denote the weak convergence in $S$ and $J_1$-topology, respectively, and
impose the following assumptions:
\begin{longlist}[(A3)]
\item[(A1)]  $S_n(t)/a_n\stackrel{S}{\Longrightarrow} S(t)$
on $D[0, 1]$. For more information about the weak convergence in the
$S$-topology, see Jakobowski (\cite{j96} and \cite{j97}),


\item[(A2)]  $g(\cdot)$ is a H\"{o}lder continuous function on
$\mathbb{R}$, that is,
$|g(x)-g(y)|\leq C|x-y|^{\nu}$ for all $x, y\in(-\infty, \infty)$,
where $\nu=1$, when $\varepsilon_1$ has infinite variance with tail index
$\alpha<2$ and $\nu>1$ when $\alpha=2$ or $\varepsilon_1$ has finite variance.

\item[(A3)]\mbox{}
\begin{enumerate}[(ii)]
\item[(i)] $ \sum_{j=1}^{\infty}
\Vert \sum_{i=j}^{\infty} F_i(x|\mathcal{F}_0)-F_i(x|\mathcal
{F}_0^*) \Vert^2<\infty,$ or

\item[(ii)] $
\sum_{i=1}^{\infty} \Vert  F_i(x|\mathcal
{F}_0)-F_i(x|\mathcal{F}_0^*) \Vert <\infty$ and
$
\sum_{i=m}^{\infty} \Vert  F_i(x|\mathcal
{F}_0)-F_i(x|\mathcal{F}_0^*) \Vert =\mathrm{O}[(\log{m})^{-a}],\break a >3/2,$
when $\sum_{i=1}^n|\varepsilon_i|/a_n=\mathrm{O}_{p}(1),$
\end{enumerate}

\item[(A4)]  $ \sum_{j=1}^{\infty}\sup_{x}
\Vert \sum_{i=j}^{\infty} F_i^{(l)}(x|\mathcal
{F}_0)-F_i^{(l)}(x|\mathcal{F}_0^*) \Vert^2<\infty, l=0,
1, F_i^{(0)}(x|\mathcal{F}_0)=F_i(x|\mathcal{F}_0).$
\end{longlist}

\begin{theorem}\label{teo2.1}
Suppose that conditions \textup{(A1)--(A3)} hold, then for any
$x\in\mathbb{R}$,
%
\begin{eqnarray}
\label{2.1}\frac{1}{\sqrt{n}}\sum_{i=1}^{n}g(S_{i-1}/a_n)
\bigl[I(\varepsilon_i<x)-F(x) \bigr]\stackrel{\mathcal {L}} {
\longrightarrow} \int_0^1 g \bigl(S(t-) \bigr)
\, \mathrm{d} W(t, x).
\end{eqnarray}
In addition, if \textup{(A3)} is replaced by \textup{(A4)},
then for any constant $A>0$,
%
\begin{eqnarray}
\label{2.2} \frac{1}{\sqrt{n}}\sum_{i=1}^{n}g(S_{i-1}/a_n)
\bigl[I(\varepsilon_i<x)-F(x) \bigr] \stackrel{w} {\Longrightarrow}
\int_0^1 g \bigl(S(t-) \bigr) \,\mathrm{d} W(t,
x), \qquad \mbox{on } D[-A, A],\qquad
\end{eqnarray}
\end{theorem}

\begin{theorem}\label{teo2.2}
Suppose that $\beta=1$ in model (\ref{1.1}) and
conditions \textup{(A1)--(A3)} in Theorem \ref{teo2.1} hold, then
%
\begin{eqnarray}
\label{2.4a}\frac{1}{\sqrt{n}}\alpha_n(x)\stackrel {\mathcal {L}} {
\longrightarrow} \int_0^1 g \bigl(S(t-) \bigr)
\, \mathrm{d} W(t, x).
\end{eqnarray}
In addition, if $a_n(\hat{\beta}-1)=\mathrm{o}_p(1)$, then
%
\begin{eqnarray}
\label{2.4b} \frac{1}{\sqrt{n}}\hat\alpha_n(x)&=&\frac
{1}{\sqrt{n}}
\alpha_n(x)
\nonumber
\\[-8pt]
\\[-8pt]
&&{}+\frac{1}{\sqrt{n}} \sum_{i=1}^{n}g(X_{i-1}/a_n)
\bigl[F \bigl(x+(\hat{\beta}-\beta )X_{i-1} \bigr)-F(x) \bigr]+
\mathrm{o}_p(1).
\nonumber
\end{eqnarray}
\end{theorem}

Let $\hat{\beta}$ be the $\tau$-quantile estimate of $\beta$
when $\varepsilon_1$ has infinite variance with tail index $\alpha<2$, that
is, $ \hat\beta=\operatorname{arg\,min}_{\beta}
\sum_{i=1}^{n}\rho_{\tau}(X_i$
$-\beta X_{i-1}-F^{-1}(\tau)),$ where $\rho_\tau(y)=y(\tau-I(y\leq
0))$. When $\beta=1$, using the argument of Theorem 4 in Knight (see
also Chan and Zhang \cite{cz09a}), we have
\[
a_n\sqrt{n}(\hat{\beta}-\beta)=\frac{({1}/{\sqrt{n}})\sum_{t=1}^{n}(X_{t-1}/a_n)  (
\tau-I (\varepsilon_t\leq F^{-1}(\tau) ) )}{
{({1}/{ n})\sum_{t=1}^{n}f_{t}
(F^{-1}(\tau)|\mathcal{F}_{t-1} )
(X_{t-1}^2/a_n^2
)}}+
\mathrm{o}_p(1).
\]
By virtue of Theorem \ref{teo2.2} and this expression, the following corollary
concerning the quantile estimate is immediate.

\begin{corollary}\label{cor2.1}
Under conditions Theorem \ref{teo2.2},
if $\mathrm{E}|f_1(F^{-1}(\tau))|\mathcal{F}_0|^{p}<\infty$ for
some $p>1$ and $f(F^{-1}(\tau))>0$, then
\[
a_n\sqrt{n}(\hat{\beta}-\beta)\stackrel{\mathcal{L}} {
\longrightarrow }-\frac{1}{ f (F^{-1}(\tau) )} \frac{\int_{0}^{1}S(t-)
\,\mathrm{d}W (t, F^{-1}(\tau) )}{\int_{0}^{1}S^2(t)
\,\mathrm{d}t}.
\]
\end{corollary}

\begin{theorem}\label{teo2.3}
In addition to the conditions of Theorem \ref{teo2.2}, if \textup{(A4)}
holds, then for any constant $A>0$,
%
\begin{eqnarray}
\label{2.5a} \sup_{x\in[-A, A]}\frac{1}{\sqrt {n}}\alpha_n(x)
\stackrel{\mathcal{L}} {\longrightarrow} \sup_{x\in[-A, A]}\int
_0^1 g \bigl(S(t) \bigr) \,\mathrm{d} W(t, x).
\end{eqnarray}
For $\hat\alpha_n(x)$, we have:
\begin{enumerate}[(b)]
\item[(a)] if $\hat{\beta}$ is the $\tau$-quantile estimate of
$\beta$ and $f(F^{-1}(\tau))>0$, then
%
\begin{eqnarray}
\label{2.5b} &&\sup_{x\in[-A, A]}\frac{1}{\sqrt {n}}\hat
\alpha_n(x)\nonumber\\
&&\quad\stackrel{\mathcal{L}} {\longrightarrow}
\sup_{x\in[-A, A]} \biggl[ \biggl(-\frac{f(x)}{ f(F^{-1}(\tau))} \biggr) \biggl(
\frac{\int_{0}^{1}S(t-)\, \mathrm{d}W(t,
F^{-1}(\tau)) }{\int_{0}^1 S^2(t)\, \mathrm{d}t} \biggr)\\
&&\hspace*{52pt}\qquad\times\int_{0}^{1}g
\bigl(S(t) \bigr)S(t)\, \mathrm{d}t+\int_0^1 g
\bigl(S(t-) \bigr) \,\mathrm{d} W(t, x) \biggr].
\nonumber
\end{eqnarray}

\item[(b)] if $\hat{\beta}$ is the LSE of $\beta$ and $(S_n(t),
\sum_{i=1}^{n}\varepsilon_i^2/a_n^2)\stackrel
{\mathit{f.d.d.}}{\longrightarrow}(S(t), S^2)$,
then
\begin{enumerate}[(ii)]
\item[(i)] if $a_n=n^\vartheta l(n)$ for some $1/2<\vartheta<1$,
%
\begin{eqnarray}
\label{2.6a}&& \sup_{x\in[-A, A]}\frac{1}{a_n}\hat \alpha_n(x)\nonumber\\
&&\quad
\stackrel{\mathcal{L}} {\longrightarrow}\sup_{x\in[-A,
A]}f(x)\int
_{0}^{1}S(t-)\, \mathrm{d}S(t)
\\
&&\hspace*{8pt}\qquad{}\times \int_{0}^{1}g \bigl(S(t) \bigr)S(t)
\, \mathrm{d}t\Big/\int_{0}^1 S^2(t) \,
\mathrm{d}t.
\nonumber\
\end{eqnarray}

\item[(ii)] If $a_n=\sqrt{n}$,
\end{enumerate}
%
\begin{eqnarray}
\label{2.6b} \hspace*{-30pt}\sup_{x\in[-A, A]}\frac{1}{\sqrt {n}}\hat
\alpha_n(x)&\stackrel{\mathcal{L}} {\longrightarrow}&
\sup_{x\in[-A,
A]} \biggl[f(x)\int_{0}^{1}S(t-)
\,\mathrm{d}S(t) \int_{0}^{1}g \bigl(S(t)
\bigr)S(t) \,\mathrm{d}t\Big/\int_{0}^1
S^2(t) \, \mathrm{d}t\hspace*{30pt}
\nonumber
\\[-8pt]
\\[-8pt]
&&\hspace*{37pt} {}+ \int_0^1 g \bigl(S(t-)
\bigr) \,\mathrm{d} W(t, x) \biggr].
\nonumber
\end{eqnarray}
\end{enumerate}
\end{theorem}

\begin{remark}\label{rem2.1}
By Voln\'{y} \cite{v93}, condition (i) in (A3) is a
necessary and sufficient condition for $I(\varepsilon_i\leq x)$
enjoying a martingale decomposition,
that is, there exist a martingale difference $\zeta_i(x)$ with respect
to $\mathcal{F}_i$ and a sequence
$\xi_i(x)\in L^2, i\in\mathcal{Z}$ such that
%
\begin{eqnarray}
\label{2.6'} I(\varepsilon_i\leq x)-F(x)=
\zeta_i(x)+\xi_{i}(x)-\xi_{i+1}(x).
\end{eqnarray}
\end{remark}

\begin{remark}\label{rem2.2}
From Theorems \ref{teo2.2} and \ref{teo2.3}, we see that the limit
distribution of the test statistics based on $\alpha_n(x)$ and $\hat
{\alpha}_n(x)$ are very different in the unit root case. As a result,
using a residual marked empirical process ($\hat{\alpha}_n(x)$) to
test the goodness-of-fit of
nonstationary processes will be very different from using the marked
empirical process ($\alpha_n(x)$).
\end{remark}

To illustrate the usefulness of these theorems, consider the following
examples, which characterize the limit
distributions of the marked empirical process $\alpha$ under various
situations.

\begin{example}\label{exe2.1}
Let $\{\varepsilon_i\}$ in model (\ref{1.1}) be the
generalized autoregressive conditional heteroscedasticity
($\operatorname{GARCH}(1, 1)$) process
\[
\varepsilon_i=\sigma_i\eta_i,\qquad
\sigma_i^2=\omega+a\sigma_{i-1}^2+b
\varepsilon_{i-1}^2,
\]
where $\omega, a, b>0$,
$\{\eta_i\}$ is an i.i.d. symmetric random sequence with $\mathrm
{E}[\log(a+b\eta_1^2)]<0$ and $\mathrm{E}(a+b\eta_1^2)^r<\infty$
for some $r>0.$
If there exists a positive constant $C_0$ such that the density
$f_{\eta}(\cdot)$ of $\eta_1$ satisfies $\sup_x f_{\eta}(x)<C_0,$
according to Kesten \cite{k73} (see also Lemma A.1 in Chan and Zhang \cite{cz10}),
there exists an $\alpha>0$ such that $\mathrm{E}(a+b\eta_1^2)^{\alpha/2}=1$ and
there exists a constant $c_0$ such that\vadjust{\goodbreak} $\lim_{x\rightarrow\infty
}x^{-\alpha}P(|\varepsilon_1|>x)=c_0.$ Let $\eta_0'$ be a
independent copy of $\eta_0$, then
there exist a constant $C$ such that
\begin{eqnarray*}
\sum_{j=1}^{\infty} \Bigg\Vert \sum
_{i=j}^{\infty} F_i(x|\mathcal{F}_0)-F_i
\bigl(x|\mathcal{F}_0^* \bigr) \Bigg\Vert^2 &\leq& \sum
_{j=1}^{\infty} \Biggl[\sum
_{i=j}^{\infty
}\big\Vert F_i(x|
\mathcal{F}_0)-F_i \bigl(x|\mathcal{F}_0^*
\bigr)\big\Vert \Biggr]^2
\nonumber
\\
&\leq& \sum_{j=1}^{\infty} \Biggl[\sum
_{i=j}^{\infty}2\Bigg\Vert \min \Biggl\{1, \frac{bC_0}{\sqrt{\omega}}
\prod_{k=1}^{i} \bigl(a+b
\eta_k^2 \bigr)\big|\eta_0^2-
\eta_0^{\prime 2}\big|\sigma_0^2 \Biggr\}
\Bigg\Vert \Biggr]^2
\nonumber
\\
&\leq& C\sum_{j=1}^{\infty} \Biggl[\sum
_{i=j}^{\infty}\Bigg\Vert \Biggl(\prod
_{k=1}^{i} \bigl(a+b\eta_k^2
\bigr) \eta_0^2\sigma_0^2
\Biggr)^{\min\{1, {\alpha/4}\}}\Bigg\Vert \Biggr]^2<\infty,
\nonumber
\end{eqnarray*}
where the last inequality follows since $\mathrm{E}(a+b\eta_1^2)^{\alpha/4}<1.$
Thus, (\ref{2.1}) of Theorem \ref{teo2.1} and Theorem \ref{teo2.2} hold (see also
Theorem 2.1 of Chan and Zhang \cite{cz10}) with
\[
a_n=\cases{ n^{1/\alpha},& \quad if $ 0<\alpha<2$,
\cr
\sqrt{ n
\log n }, & \quad if $ \alpha=2$,
\cr
\sqrt{n}, & \quad if $\alpha>2$. }
\]
Further, if $f_{\eta}(\cdot)$ has derivative $f'_{\eta}(\cdot)$ and
$\sup_x f'_{\eta}(x)<C_0,$ then
\begin{eqnarray*}
\sum_{j=1}^{\infty}\sup_{x} \Bigg\Vert
\sum_{i=j}^{\infty} f_i(x|
\mathcal{F}_0)-f_i \bigl(x|\mathcal{F}_0^*
\bigr) \Bigg\Vert^2 &\leq& C\sum_{j=1}^{\infty}
\Biggl[\sum_{i=j}^{\infty}\rho^i
\Biggr]^2<\infty.
\end{eqnarray*}
From Theorem \ref{teo2.1} it follows that for any constant $A>0$,
\begin{eqnarray*}
\frac{1}{\sqrt{n}}\sum_{i=1}^{n}g(S_{i-1}/a_n)
\bigl[I(\varepsilon_i<x)-F(x) \bigr] \stackrel{w} {\Longrightarrow}
\int_0^1 g \bigl(S(t-) \bigr) \,\mathrm{d} W(t,
x), \qquad \mbox{on } D[-A, A],
\end{eqnarray*}
where $S(t)$ is an $\alpha$-stable process when $\alpha<2$ and a
Gaussian process when $\alpha\geq2$
and
$W(t, x)$ is given as in Theorem \ref{teo2.1}. As a result, when $\beta=1$ in
model (\ref{1.1}),
\[
\sup_{x\in[-A, A]} \frac{1}{\sqrt{n}}\alpha_n(x)\stackrel {
\mathcal{L}} {\longrightarrow} \sup_{x\in[-A, A]}\int_0^1
g \bigl(S(t) \bigr) \,\mathrm{d} W(t, x).
\]
\end{example}

\begin{example}\label{exe2.2}
Let $\{\varepsilon_i\}$ in model (\ref{1.1}) be an
infinite-variance linear moving average process $\varepsilon_i=\sum_{j=1}^{\infty}c_j\eta_{i-j}$, where $ \{\eta_i\} $ is an i.i.d.
sequence with bounded density and heavy tail index $0<\alpha\leq2$,
that is, when $\alpha<2$,
$nP(|\eta_1|>a_n x)\rightarrow x^{-\alpha}$ for any $x>0$ and when
$\alpha=2$, $na_n^{-2}\mathrm{E} (\eta_1^{2}I(|\eta_1|\leq
a_n))\rightarrow1$. Since
%
\begin{eqnarray}
\label{2.7} &&\sum_{j=1}^{\infty} \Bigg\Vert \sum
_{i=j}\bigl[ F_i(x|
\mathcal{F}_0)-F_i \bigl(x|\mathcal{F}_0^*
\bigr)\bigr] \Bigg\Vert^2
\nonumber
\\
&&\quad \leq\sum_{j=1}^{\infty} \Biggl[\sum
_{i=j}^{\infty}\big\Vert F_i(x| \mathcal
{F}_0)-F_i \bigl(x|\mathcal{F}_0^* \bigr)
\big\Vert \Biggr]^2
\\
&&\quad \leq\sum_{j=1}^{\infty} \Biggl[\sum
_{i=j}^{\infty}C\big\Vert \min \bigl(\big|c_i
\bigl( \eta_0-\eta'_0 \bigr)\big|, 1 \bigr)
\big\Vert \Biggr]^2 \leq C'\sum
_{j=1}^{\infty} \Biggl[\sum_{i=j}^{\infty}|c_i|^{\alpha
/2}
\Biggr]^2,
\nonumber
\end{eqnarray}
if $\sum_{j=1}^{\infty}\sum_{i=j}^{\infty}|c_i|^{\alpha/2}<\infty$
 (i.e., $\sum_{i=1}^{\infty} i|c_i|^{\alpha/2}<\infty)$, condition
(i) of (A3) holds.
Further, by Chan and Zhang \cite{cz09b} (see also Avram and Taqqu \cite{at92}),
we also have condition (A1) for $\{\varepsilon_i\}$.
Thus, if $\sum_{j=1}^{\infty}\sum_{i=j}^{\infty}|c_i|^{\alpha/2
}<\infty$ and condition (A2) holds, then
for the unit-root model (1.1), we have for any $x\in\mathbb{R},$
%
\begin{eqnarray}
\label{2.8}\frac{1}{\sqrt{n}}\alpha_n(x)\stackrel {\mathcal {L}} {
\longrightarrow} \int_0^1 g
\bigl(Z_\alpha(t-) \bigr)\, \mathrm{d} W(t, x),
\end{eqnarray}
where $Z_\alpha(t)$ is a stable process with index $\alpha$. In
particular, when $c_j=j^{-\theta}$ and $\theta>3/\alpha$, then under
condition (A2), conclusion (\ref{2.8}) holds.

On the other hand, if $0<\alpha<1$, since $\sum_{i=1}^{n}|\varepsilon_i|/a_n=\mathrm{O}_p(1),$ using (ii) of condition (A3), we can relax the condition
from $\theta>3/\alpha$ to $\theta>2/\alpha.$ This observation sheds
light on the important subtlety of the roles of $\theta$ and $\alpha$
for an infinite variance moving average process.
\end{example}

\begin{example}\label{exe2.3}
When $ \{\eta_i\} $ in Example \ref{exe2.2} has finite variance
and bounded density $f_\eta(x)$ and $c_j=j^{-\theta}l(j),$ as $
j\rightarrow\infty,$ for some slowly varying function, similar to
(\ref{2.7}), we have that under $\theta>1$, Theorem \ref{teo2.2} holds.
Further, if $f'_\eta(x)$ exists and $\sup_x |f'_\eta(x)|\leq C_0$
for some $C_0>0$ and $\theta>3/2$, then
\begin{eqnarray*}
&&\sum_{j=1}^{\infty} \Bigg\Vert \sum
_{i=j}^{\infty} f_i(x|\mathcal{F}_0)-f_i
\bigl(x|\mathcal{F}_0^* \bigr) \Bigg\Vert^2 \\
&&\quad \leq\sum
_{j=1}^{\infty} \Biggl[\sum
_{i=j}^{\infty}C\big\Vert \min \bigl(\big|c_i \bigl(
\eta_0-\eta'_0 \bigr)\big|, 1 \bigr)\big\Vert
\Biggr]^2
\\
&&\quad \leq C'\sum_{j=1}^{\infty}
\Biggl[\sum_{i=j}^{\infty}i^{-\theta
}l(i)
\Biggr]^2 \leq C''\sum
_{j=1}^{\infty}i^{-2\theta+2}l'(i)<\infty,
\end{eqnarray*}
where $l'(x)$ is a slowly varying function, it follows that condition
(A4) holds. Thus, if condition (A2) holds,
\begin{eqnarray*}
\frac{1}{\sqrt {n}}\sum_{i=1}^{n}g(S_{i-1}/
\sqrt{n}) \bigl[I(\varepsilon_i<x)-F(x) \bigr] \stackrel{\mathcal{L}}
{ \longrightarrow} \int_0^1 g \bigl(S(t)
\bigr) \,\mathrm{d} W(t, x),\qquad \mbox{on } D[-A, A]
\end{eqnarray*}
and when $\beta=1$ in model (\ref{1.1}),
\[
\sup_{x\in[-A, A]} \frac{1}{\sqrt{n}}\sum_{i=1}^{n}g(X_{i-1}/
\sqrt{n}) \bigl[I(\varepsilon_i<x)-F(x) \bigr]\stackrel{\mathcal {L}}
{ \longrightarrow} \sup_{x\in[-A, A]}\int_0^1
g \bigl(S(t) \bigr)\, \mathrm{d} W(t, x),
\]
where $S(t)$ is a Gaussian process.
\end{example}

\section{Long-memory error processes}\label{sec3}

In this section, we study the marked empirical process
\[
\alpha_n(x)=\sum_{i=1}^{n}g(X_{i-1}/a_n)
\bigl(I(\varepsilon_i<x)-F(x) \bigr),
\]
when $\{X_t\}$ is a unit root process given by (\ref{1.1}) with $\beta
=1$ and $\varepsilon_t$ being a long-memory process.
Long-memory processes have been widely applied in finance and
econometrics, see, for example, Baillie \cite{b96} and Teyssi\'{e}re and
Kirman \cite{tk07}.
Specifically, let $c_j=j^{-\theta}l(j), l(\cdot)$ be a slowly
varying function with $| l(m+n)/l(n) -1 |\leq m/n$ for $1\leq m\leq n$ and consider the linear moving average process
$\varepsilon_i=\sum_{j=1}^{\infty}c_j\eta_{i-j}$ defined in Example~\ref{exe2.3} with
$\sum_{j=1}^{\infty}c_j^2<\infty$ and $\sum_{j=1}^{\infty
}|c_j|=\infty$.

The essential idea in studying the weak convergence of
$\alpha_n(x)$ when $\{\varepsilon_i\}$ is short-memory is the
martingale approximation.
This transforms the weak convergence of $\alpha_n(x)$ into those of a
martingale stochastic integral $\sum_{i=1}^{n}g(X_{i-1}/a_n)\xi_i(x)$.
When $\{\varepsilon_i\}$ is long-memory, this method does not work and
the issue of the weak convergence
of $ \alpha_n(x)$ becomes much more challenging. Fortunately, to
circumvent this difficulty, the ideas of Ho and Hsing \cite{hh96} and Wu
\cite{w03} become relevant.

Let $f(x)$ be the density of $\eta_1$ and $f^{(l)}(x)$ be its $l$th
derivatives and let $f(x)=f^{(0)}(x)$. We have the following results.

\begin{theorem}\label{teo3.1}
Suppose that $\beta=1$ in model (\ref{1.1}) and \textup{(i)}
$E|\eta_1|^{\nu}<\infty$ for some $\nu>\max\{4, 1/(1-\theta)\};$ \textup{(ii)} $g(\cdot)$ is Lipschitz
on $\mathbb{\mathbb{R}};$ \textup{(iii)} $\sum_{l=0}^{p}\int_{R}|f^{(l)}(x)|^2 \,\mathrm{d}x<\infty$, then if any one of the following
three conditions holds:
\begin{enumerate}[(b)]
\item[(a)] $p=4$ and $\theta\in(1/2, 3/4)\cup(5/6, 1);$

\item[(b)] $p=5$ and $\theta\in(1/2, 3/4)\cup(3/4, 1)$ or $\theta=3/4$ and
$\sum_{i=1}^{\infty}l^4(i)/ i<\infty;$

\item[(c)] $p=6$ and $\theta\in(1/2, 1),$
\end{enumerate}
we have
%
\begin{eqnarray}
\label{3.1}\sup_{x}\frac{1}{ a_n}\alpha_n(x)\stackrel
{\mathcal{L}} {\longrightarrow} \sup_{x}\int_0^1
g \bigl(Z_{\theta}(t) \bigr) \,\mathrm{d} Z_{\theta}(t, x).
\end{eqnarray}
If in addition $n(\hat\beta-1)=\mathrm{O}_p(1)$, then
%
\begin{eqnarray}
\label{3.2} &&\sup_{x}\frac{1}{ a_n} \Biggl[\hat
\alpha_n(x)-f(x) (\hat \beta-\beta)\sum_{i=1}^{n}g(X_{i-1}/a_n)X_{i-1}
\Biggr]
\nonumber
\\[-8pt]
\\[-8pt]
&&\quad \stackrel{\mathcal{L}} {\longrightarrow} \sup_{x}f(x)\int
_0^1 g \bigl(Z_{\theta}(t) \bigr)\,
\mathrm{d} Z_{\theta}(t),
\nonumber
\end{eqnarray}
where $a_n=n^{3/2-\theta}l(n)$ and $Z_{\theta}(t)=\int_{-\infty
}^{t}\int_{0}^{t}[\max(v-u, 0)]^{-\theta} \,\mathrm{d}v \,\mathrm{d}B(u)$, $B(u)$ is a
standard Brownian motion.
\end{theorem}

\begin{remark}\label{rem3.1}
When $g(\cdot)\equiv1$, then Theorem \ref{teo3.1} reduces to
the case of Chan and Ling~\cite{cl08}.
\end{remark}

\section{Proofs}\label{sec4}

To prove the main results, we need the following lemmas. The first one
is due to Lemma 4 of Wu~\cite{w03}.

\begin{lemma}\label{lem4.1}
Let $H\in\mathcal{C}^1, $ the space of functions with
continuous first-order derivatives and \mbox{$a>0$}. Then
%
\begin{eqnarray}
\sup_{t\leq s\leq t+a}H^2(s)\leq{2\over a}\int_{t}^{t+a}H^2(u)
\,\mathrm{d}u + 2a\int_{t}^{t+a}H'^2(u)
\,\mathrm{d}u
\end{eqnarray}
and
%
\begin{eqnarray}
\sup_{t\in\mathbb{R}}H^2(s)\leq2\int_{\mathbb{R}}H^2(u)
\,\mathrm{d}u + 2\int_{\mathbb{R}}H'^2(u)
\,\mathrm{d}u,
\end{eqnarray}
where $H'$ is the derivative of $H$.
\end{lemma}

\begin{lemma}\label{lem4.2}
If $\{\varepsilon_i\}$ is short-memory, then under the
conditions \textup{(A1)}, \textup{(A2)} and \textup{(A3)}, there exists a martingale difference
sequence $\zeta_i(x)$ with respect to $\mathcal{F}_i$ such that for
any $\delta>0$,
\[
\lim_{n\rightarrow\infty} P \Biggl\{ \Bigg| \frac{1}{\sqrt{n}}\sum
_{i=1}^{n}g(S_{i-1}/a_n)
\bigl(I(\varepsilon_i<x)-F(x) \bigr)-\frac{1}{\sqrt{n}}\sum
_{i=1}^{n} g(S_{i-1}/a_n)
\zeta_i(x) \Bigg|>\delta \Biggr\}=0.
\]
\end{lemma}

\begin{pf}
When (i) of (A3) holds, then by Voln\'{y} \cite{v93}, there
exist a random sequence\linebreak[4]  $\xi_i(x)=\sum_{j=-\infty}^{-1}\sum_{l=0}^{\infty}\mathcal{P}_{i+j}I(\varepsilon_{i+l}\leq x)\in L^2$
and a
martingale difference sequence\linebreak[4]  $\zeta_i(x)= \sum_{j=i}^{\infty
}{\mathcal{P}}_{i}I(\varepsilon_j\leq x)$ such that
$I(\varepsilon_i<x)-F(x)=\zeta_i(x)+\xi_i(x)-\xi_{i+1}(x).$
This gives that
%
\begin{eqnarray}
\label{4.3} &&\frac{1}{\sqrt{n}}\sum_{i=1}^{n}g(S_{i-1}/
a_n) \bigl(I(\varepsilon_i<x)-F(x) \bigr)-
\frac{1}{\sqrt{n}}\sum_{i=1}^{n}
g({S_{i-1}/a_n})\zeta_i(x)
\nonumber
\\[-8pt]
\\[-8pt]
&&\quad = \frac{1}{\sqrt{n}}\sum_{i=1}^{n-1}
\xi_{i+1}(x) \bigl[g({S_{i}/ a_n})-g({S_{i-1}/
a_n}) \bigr]-\frac{1}{\sqrt{n}}g(S_{n-1}/a_n)
\xi_{n+1}(x) =:I_1+I_2.\qquad
\nonumber
\end{eqnarray}
Since for any $\delta>0$,
%
\begin{eqnarray}
\label{4.4} P \Bigl\{\sup_{2\leq i\leq n+1}\big|\xi_{i}(x)\big|>\delta\sqrt {n}
\Bigr\} 
&\leq&\sum_{i=2}^{n+1} (
\sqrt{n}\delta)^{-2}\mathrm{E} \bigl[\xi_{1}^2(x)I
\bigl(\big|\xi_{1}(x)\big|>\delta\sqrt{n} \bigr) \bigr]\rightarrow0,
\end{eqnarray}
it follows that $|\xi_{n+1}(x)|/\sqrt{n}=\mathrm{o}_p(1)$. On the other hand,
by (A1) and (A2), we have $g(S_{n-1}/a_n)=\mathrm{O}_p(1).$
Thus, $I_2=\mathrm{o}_p(1).$ It suffices to show that $I_1=\mathrm{o}_p(1)$. When $\{
\varepsilon_i\}$ has infinite variance with tail index
$\alpha<2$, the result $I_1=\mathrm{o}_p(1)$ follows
along exactly the lines of argument of Lemma 2 of Knight \cite{k91}. We
therefore only give the proof for the finite variance case in detail.

When $\{\varepsilon_i\}$ has finite variance or has infinite variance with tail index $\alpha=2$, since $\sum_{i=0}^{\infty}\theta_2(i)<\infty$, it follows from Theorem 1 of Wu
\cite{w07} that
$\mathrm{E}(\sum_{i=1}^{n}\varepsilon_i)^2= C_1n.$ Thus,
$a_n=C_2\sqrt{n}.$ By (A2) and (i) of (A3) for any $\epsilon>0$, we have
%
\begin{eqnarray}
P(|I_1|>\epsilon)&\leq& {C\over \sqrt{n}}\sum_{i=1}^{n-1} E \bigl|\xi_{i+1}(x)(\varepsilon_i/a_n)^{\nu}
I\bigl(|\varepsilon_i|\leq \delta a_n\bigr)\bigr|+P\Bigl(\sup_{1\leq i\leq n}|\varepsilon_1|>\delta a_n\Bigr)
\nonumber
\\
&\leq& {C'\delta^{\nu-1}\over n}\sum_{i=1}^{n-1}
\bigl\{ \bigl[\mathrm{E}\xi_{i+1}^2(x) \bigr]^{1/2}
\bigl[\mathrm{E}\varepsilon_{i}^2 \bigr]^{1/2}
\bigr\} \leq{C''\delta^{\nu-1}\over n}\sum
_{i=2}^{n} \bigl\{ \bigl[\mathrm{E}
\xi_{i}^2(x) \bigr]^{1/2}\bigr\}
\\
&\leq& {C'''\delta^{\nu-1}\over n}\sum
_{i=2}^{n} \Biggl\{\sum_{j=-\infty}^{i-1}\mathrm{E} \Biggl[\sum
_{l=0}^{\infty
} \bigl[F_{i+l}(x|
\mathcal{F}_j)-F_{i+l} \bigl(x|\mathcal{F}_{j}^{*}
\bigr) \bigr] \Biggr]^2 \Biggr\}^{1/2}=\mathrm{o}(1)
\nonumber
\end{eqnarray}
by taking $\delta\rightarrow 0$.
This gives that $I_1=\mathrm{o}_p(1)$ and therefore Lemma \ref{lem4.2} holds when (i) of
(A3) is true.

When (ii) of (A3) holds, by Corollary 1 of Wu \cite{w07}, we have
%
\begin{eqnarray}
\label{4.6}\sup_{0\leq t\leq1} \Bigg|W_n(t, x)-\sum
_{i=1}^{[nt]}\zeta_i(x) \Bigg|=\mathrm{o}
\bigl(n^{1/2} \bigr),\qquad \mbox{a.s.}
\end{eqnarray}
where $W_n(t, x)=\sum_{i=1}^{[nt]}(I(\varepsilon_i\leq x)-F(x))$.
Combining this with (A2) gives
%
\begin{eqnarray}
&& \Bigg|\frac{1}{\sqrt{n}}\sum_{i=1}^{n}g(S_{i-1}/a_n)
\bigl(I(\varepsilon_i<x)-F(x) \bigr)-\frac{1}{\sqrt {n}}\sum
_{i=1}^{n} g(S_{i-1}/a_n)
\zeta_i(x) \Bigg|
\nonumber
\\
&&\quad =\frac{1}{\sqrt{n}} \Bigg|\sum_{i=1}^{n}
\Biggl(W_n(i/n, x)-\sum_{i=1}^{i}
\zeta_i(x) \Biggr) \bigl[g(S_{i}/a_n)-g(S_{i-1}/a_n)
\bigr] \Bigg|
\\
&&\quad \leq C \Biggl(\sup_{0\leq t\leq1}\frac{1}{\sqrt{n}} \Bigg|W_n(t,
x)- \sum_{i=1}^{[nt]}\zeta_i(x)
\Bigg|\Biggr) \Biggl(\sum_{i=1}^{n}
\frac{|
\varepsilon_i|}{ a_n} \Biggr)=\mathrm{o}_p(1).
\nonumber
\end{eqnarray}
This completes the proof of Lemma \ref{lem4.2}.
\end{pf}

\begin{lemma}\label{lem4.3}
If $\{\varepsilon_i\}$ is short-memory, then under the
conditions \textup{(A1)}, \textup{(A2)} and \textup{(A4)},
for any constant $A>0$,
\[
\sup_{x\in[-A, A]} \Bigg| \frac{1}{\sqrt{n}}\sum_{i=1}^{n}g(S_{i-1}/a_n)
\bigl(I(\varepsilon_i<x)-F(x) \bigr)-\frac{1}{\sqrt{n}}\sum
_{i=1}^{n} g(S_{i-1}/a_n)
\zeta_i(x) \Bigg|
\]
converges to zero in probability,
where $\zeta_i(x)$ is defined as in Lemma \ref{lem4.2}.
\end{lemma}

\begin{pf}
From the proof of Lemma \ref{lem4.2} for case of (i) of (A3), it
suffices to show
%
\begin{eqnarray}
\label{4.7}\frac{1}{ n}\sum_{i=1}^{n}
\mathrm{E}\sup_{x\in[-A,
A]}\xi_i^2(x)=\mathrm{O}(1).
\end{eqnarray}
Since $\xi_i(x)=\sum_{j=-\infty}^{i-1}\sum_{l=0}^{\infty}\mathcal
{P}_{j}I(\varepsilon_{i+l}\leq x)\in C^1$, it follows from
Lemma \ref{lem4.1} and Fubini's theorem that
\begin{eqnarray*}
\mathrm{E}\sup_{x\in[-A, A]}\xi_i^2(x)&\leq&
\frac{2}{A} \int_{-A}^{A}\mathrm{E}
\xi_i^2(u) \,\mathrm{d}u + 2A\int_{-A}^{A}
\mathrm{E} \xi_i'^2(u) \,\mathrm{d}u
\nonumber
\\
&\leq& \frac{2}{A}\int_{-A}^{A}\sum
_{j=-\infty}^{i-1} \mathrm{E} \Biggl[\sum
_{l=0}^{\infty} \bigl[F_{i+l}(u|\mathcal
{F}_j)-F_{i+l} \bigl(u|\mathcal{F}_{j}^{*}
\bigr) \bigr] \Biggr]^2 \,\mathrm{d}u
\nonumber
\\
&&{}+2A\int_{-A}^{A}\sum
_{j=-\infty}^{i-1} \mathrm{E} \Biggl[\sum
_{l=0}^{\infty} \bigl[f_{i+l}(u|\mathcal
{F}_j)-f_{i+l} \bigl(u|\mathcal{F}_{j}^{*}
\bigr) \bigr] \Biggr]^2 \,\mathrm{d}u
\nonumber
\\
&\leq&2\sum_{j=-\infty}^{i-1} \sup_{u}
\mathrm{E} \Biggl[\sum_{l=0}^{\infty}
\bigl[F_{i+l}(u|\mathcal {F}_j)-F_{i+l} \bigl(u|
\mathcal{F}_{j}^{*} \bigr) \bigr] \Biggr]^2
\nonumber
\\
&&{}+4A^2\sum_{j=-\infty}^{i-1}
\max_{u} \mathrm{E} \Biggl[\sum_{l=0}^{\infty}
\bigl[f_{i+l}(u|\mathcal {F}_j)-f_{i+l} \bigl(u|
\mathcal{F}_{j}^{*} \bigr) \bigr] \Biggr]^2.
\nonumber
\end{eqnarray*}
Thus, by (A4), we have (\ref{4.7}) as desired.
\end{pf}

\begin{lemma}\label{lem4.4}
Let $\widetilde{W}_n(t, x)=\sum_{i=1}^{[nt]}\zeta_i(x), \zeta_i(x)$ is the martingale difference defined in
Lem\-ma~\ref{lem4.2}. Then under condition \textup{(A4)},
\[
\frac{1}{\sqrt{n}}\widetilde{W}_n(t, x)\stackrel{w} {\Longrightarrow}
W(t, x)\qquad \mbox{on } D \bigl([0, 1]\times[-A, A] \bigr).
\]
\end{lemma}

\begin{pf}
Condition (A4) implies $ \sum_{j=1}^{\infty}\Vert \sum_{i=j}^{\infty} F_i(x|\mathcal{F}_0)-F_i(x|\mathcal
{F}_0^*)\Vert^2<\infty$, it follows that $E\zeta_i(x)=
E\{\sum_{i=0}^{\infty} F_i(x|\mathcal{F}_0)-F_i(x|\mathcal{F}_0^*)\}^2=\mu(x)<\infty$. Since $\{\zeta_i(x)\}$ is a martingale difference
sequence, by
Theorem 23.1 of Billingsley \cite{b68}, we have
%
\begin{eqnarray}
\label{4.8}\frac{1}{\sqrt{n}}\widetilde{ W}_n(t, x)\stackrel {
\mathcal{L}} {\longrightarrow} W(t, x).
\end{eqnarray}
%
By (\ref{4.8}) and the Cram\'{e}r--Wold's device, the
finite-dimensional convergence of $\widetilde{W}_n(t, x)$ follows.
By Theorem 6 of Bickel and Wichura \cite{bw71}, to show the tightness of $\{
\widetilde{W}_n(t, x)\}$ on $D[0, 1]\times D[-A, A]$, it suffices to
show that for any $0\leq t_1<t<t_2\leq1$ and $-A\leq x_1< x< x_2\leq A$,
%
\begin{eqnarray}
\label{4.9a} n^{-2}E \Biggl\{ \Biggl[\sum
_{i=[nt_1]+1}^{[nt]} \zeta_i(x_1,
x_2) \Biggr]^2 \Biggl[\sum_{i=[nt]+1}^{[nt_2]}
\zeta_i(x_1, x_2) \Biggr]^2 \Biggr
\}\leq(t-t_1) (t_2-t) (x_2-x_1)^2
\quad
\end{eqnarray}
and
%
\begin{eqnarray}
\label{4.9b}n^{-2} \mathrm{E} \Biggl\{ \Bigg|\sum_{i=[nt_1]+1}^{[nt_2]}
\zeta_i(x_1, x) \Bigg|^2 \Bigg|\sum
_{i=[nt_1]+1}^{[nt_2]}\zeta_i(x, x_2)
\Bigg|^2 \Biggr\} \leq C(x-x_1) (x_2-x)
(t_2-t_1)^2,
\end{eqnarray}
where $\zeta_i(x, y)=\zeta_i(y)-\zeta_i(x).$ Equations (\ref{4.9a})
and (\ref{4.9b}) follow easily by condition (A4) and noting that
$\widetilde{W}_n(t, x)$ is a martingale. Details are omitted.
\end{pf}

\begin{lemma}\label{lem4.5}
Under the conditions of Theorem \ref{teo2.1}, there exists a dense
set $Q\subset[0, 1], 0, 1\in Q$ such that
for any finite subset $\{0\leq t_1<t_2<\cdots< t_m\leq1\}\subset Q$
and for any $x,$
%
\begin{eqnarray}
\label{4.10} \bigl(S_n(t_i)/a_n, \widetilde{W}_n(t_i, x)/\sqrt{n}, 1\leq i\leq m \bigr)\stackrel{\mathcal{L}} {\longrightarrow}
\bigl(S(t_i), W(t_i, x), 1\leq i\leq m \bigr).
\end{eqnarray}
\end{lemma}

\begin{pf} Since $S_n(t)\stackrel{S}{\Longrightarrow} S(t)$, it
follows that there exists a dense set $Q'\subset[0, 1], 1\in Q'$ such that
for any finite subset $\{ t_1<t_2<\cdots< t_m\leq1\}\subset Q'$,
%
\begin{eqnarray}
\label{4.11} a_n^{-1}\bigl(S_n(t_1), S_n(t_2),
\ldots, S_n(t_m) \bigr)\stackrel {\mathrm{f.d.d.}} {
\longrightarrow} \bigl(S(t_1), S(t_2), \ldots,
S(t_m) \bigr).
\end{eqnarray}
Note that $S_n(0)=S(0)=0$. Thus, (\ref{4.11}) holds for all finite
subset of $Q=Q'\cup\{0\}.$

When $\varepsilon_i$ has infinite variance,
since $W(t, x)$ is a continuous process on $[0, 1]\times[-A, A]$, it
follows (see, e.g.,
page 112 of Billingsley \cite{b68}) that the weak convergence of Lemma~\ref{lem4.4} can
also be replaced by $C([0, 1]\times C[-A, A])$. Thus, $(S_n(t),
\widetilde{W}_n(t, x))$ is uniformly $S$-tight on $D[0, 1]\times D[0,
1]$. This implies that for any sequence $ (S_n(t), \widetilde{W}_n(t,
x)), t\in Q$, there exists a subsequence $(S_{nk}(t), \widetilde
{W}_{nk}(t, x))$ such that
\[
\mathbf{Z}_{nk}(t):= \bigl(S_{nk}(t)/a_n,
\widetilde{W}_{n}(t, x)/\sqrt{n} \bigr)\stackrel{\mathcal{L}} {\longrightarrow}
\mathbf{Z}(t),
\]
where
$ \mathbf{Z}(t)$ is a bivariate random process with marginal distributions
$S(t)$ and $W(t, x).$ Following the argument
of Theorem 3 in Resnick and Greenwood \cite{rg79}, we have that $S(t)$ and
$W(t, x)$ are independent and any convergent subsequence has the same limit.
Thus, (\ref{4.10}) holds.

When $\varepsilon_i$ has finite variance, since $\sum_{i=0}^{\infty
}\theta_2(i)<\infty$, it follows from
Corollary 3 of Dedecker and Merlev\`{e}de \cite{dm03} that
$S_n(t)/a_n\stackrel{w}{\Longrightarrow} S(t)$ for some Gaussian process
$S(t).$ Thus, by Lemma~\ref{lem4.4}, if we can show that for any finite subset
$\{t_i, 1\leq i\leq m\}\subset[0, 1]$,
%
\begin{eqnarray}
\label{4.12} \bigl(S_n(t_i)/a_n, \widetilde{W}_n(t_i, x)/\sqrt{n}, 1\leq i\leq m \bigr)\stackrel{\mathcal{L}} {\longrightarrow}
\bigl(S(t_i), W(t_i, x), 1\leq i\leq m \bigr),
\end{eqnarray}
then $(S_{n}(t)/a_n, \widetilde{W}_{n}(t, x)/\sqrt{n})\stackrel{w}{\Longrightarrow}
(S(t), W(t, x))$ on $D[0, 1]$ and (\ref{4.10}) follows. By Theorem 1
of Wu \cite{w07}, we have that
there exists martingale $E_i$ with respect to $\mathcal{F}_i$ such that
%
\begin{eqnarray}
\label{4.14} \Bigg| \bigl(S_{n}(t)/a_n, \widetilde{W}_{n}(t, x)/\sqrt{n}
\bigr)- \Biggl(\sum_{i=1}^{[nt]}E_i/a_n,
\sum_{i=1}^{[nt]}\zeta_i(x)/\sqrt{n}
\Biggr) \Bigg|=\mathrm{o}_p(1).
\end{eqnarray}
On the other hand, from the martingale central limit theorem (see
Theorem 4.1 of Hall and Heyde \cite{hh80}), it follows that
%
\begin{eqnarray}
\label{4.15}\sum_{i=1}^{[nt]}
\bigl(E_i/a_n, \zeta_i(x)/\sqrt{n} \bigr)\stackrel {w} {
\Longrightarrow} \bigl(S(t), W(t, x) \bigr).
\end{eqnarray}
Combining (\ref{4.14}) with (\ref{4.15}) yields (\ref{4.12}). This
completes the proof of Lemma \ref{lem4.5}.
\end{pf}

\begin{lemma}\label{lem4.6}
Under the conditions of Theorem \ref{teo3.1}, we have
\begin{enumerate}[(b)]
\item[(a)] $ \frac{1}{ a_n}S_n(t)=\frac{1}{ a_n}\sum_{i=1}^{[nt]}\varepsilon_i\stackrel{w}{\Longrightarrow}Z_{\theta
}(t) \mbox{ on } D[0, 1]$;

\item[(b)] $ \frac{1}{ a_n}\sum_{i=1}^{n}g(S_{i-1}/a_n)\varepsilon_i\stackrel{\mathcal
{L}}{\longrightarrow}\int_{0}^1 g(Z_\theta(t))
\,\mathrm{d}Z_\theta(t)$.
\end{enumerate}
\end{lemma}

\begin{pf}
(a) can be found in Avram and Taqqu \cite{at87}. Next, we give
the proof of (b).

With the help of strong approximation, it can be shown that
\[
{1\over a_n}\sum_{i=1}^{n}g(S_{i-1}/a_n)\epsilon_i={1\over a_n}\sum_{i=1}^{n}g(S^*_{i-1}/a_n)\epsilon^*_i +o_p(1),
\]
where $S^*_i=\sum_{j=1}^{i}\epsilon^*_j$ and $\epsilon^*_j$ is defined similarly to $\epsilon_j$ by replacing $\{\eta_i\}$ with
i.i.d normal variables $\{\eta^*_i\}$.

Since $Z_{\theta}(t)$ is a fractional Brownian motion with Hurst index
$H=3/2-\theta$, by Theorem 4 of Marcus \cite{m68}, we have
\[
P \Bigl(\limsup_{|s-t|=h\rightarrow0; 0\leq s, t\leq1}\big|Z_{\theta
}(s)-Z_{\theta}(t)\big|<2h^{H}
\log(1/h) \Bigr)=1.
\]
Thus, by (a) and $Z_{\theta}(t)$ is continuous, we have
\begin{eqnarray*}
&&\lim_{n\rightarrow\infty}P \biggl\{\limsup_{|s-t|\leq h\rightarrow 0;
0\leq s, t\leq1}\frac{1}{ a_n}\big|S^{*}_n(t)-S^{*}_n(s)\big|
\leq2h^{H}\log(1/h) \biggr\}
\nonumber
\\
&&\quad = P \Bigl(\limsup_{|s-t|\leq h\rightarrow 0; 0\leq s, t\leq
1}\big|Z_{\theta}(t)-Z_{\theta}(s)\big|
\leq2h^{H}\log(1/h) \Bigr)=1.
\nonumber
\end{eqnarray*}
This implies in probability $S^{*}_n(t)$ is H\"{o}lder continuous with an
exponent $a>H.$ This gives that, in probability, for any $p>(1/H,
\infty)$,
\[
\nu_p \bigl(S^{*}_n(t)/a_n, [0, 1] \bigr)=
\sup_{\kappa}\sum_{i=1}^m\big|S^{*}_n(t_i)-S^{*}_n(t_{i-1})\big|^{p}/a_n^p<
\infty,
\]
where the supremum is taken over all subdivisions $\kappa$
of $[0, 1]: 0=x_0<\cdots<x_m=1, m\geq1.$ Since $g(\cdot)$ is a
Lipschitz function, we have in probability
\[
\nu_p \bigl(g \bigl(S^{*}_n(t)/a_n \bigr), [0,
1] \bigr)=\sup_{\kappa}\sum_{i=1}^m\big|g
\bigl(S^{*}_n(t_i)/a_n \bigr)-g
\bigl(S^{*}_n(t_{i-1})/a_n \bigr)\big|^{p}<
\infty.
\]
By the theorem on Stieltjes integrability of Young \cite{y36} (see also
Theorem 2.4 of Mikosch and Norvai\u{s}a \cite{mn00}), we have in probability
the integral
\[
\frac{1}{ a_n}\sum_{i=1}^{n}g(S_{i-1}/a_n)
\varepsilon_i=\int_0^1 g
\bigl(S_n(t-)/a_n \bigr) \,\mathrm{d} S_n(t)/a_n
\]
exists. This implies that
%
\begin{eqnarray}
\label{4.16}&&\int_0^1 g \bigl(S_n(t-)/a_n
\bigr) \,\mathrm{d} S_n(t)/a_n
\nonumber
\\[-8pt]
\\[-8pt]
\nonumber
&&\qquad=\lim_{\delta
\rightarrow0}\sum
_{i=1}^{m}g \bigl(S_n(t_{i})
\bigr) \bigl(S_{n}(t_{i+1})-S_n(t_i)
\bigr)/a_n
\end{eqnarray}
for some
sub-division $\kappa$
of $[0, 1]\dvt 0=t_0<t_1<\cdots<t_m\leq1, m=[1/\delta]$ with
$t_{i+1}-t_i=\delta.$ By (a) and the continuous mapping theorem, we
get that for any given $m$,
%
\begin{eqnarray}
\label{4.17}\sum_{i=1}^{m}g
\bigl(S_{n}(t_i)/a_n \bigr) \bigl(S_{n}(t_{i+1})-S_n(t_i)
\bigr)/a_n&\stackrel {\mathcal{L}} {\longrightarrow}&\sum
_{i=1}^{m} g \bigl(Z_\theta(t_{i})
\bigr) \bigl(Z_\theta(t_{i+1})-Z_\theta(t_i)
\bigr)
\nonumber
\\[-8pt]
\\[-8pt]
&\stackrel{p} {\longrightarrow}&\int_0^1 g
\bigl(Z_\theta(t) \bigr) \,\mathrm{d}Z_\theta (t),
\nonumber
\end{eqnarray}
where the last equality is followed by taking $\delta\rightarrow0$
and the existence of $\int_0^1 g(Z_\theta(t)) \,\mathrm{d}Z_\theta(t).$
Combining (\ref{4.16}) and (\ref{4.17}) gives (b). The proof of Lemma
\ref{lem4.6} is completed.
\end{pf}

\begin{pf*}{Proof of Theorem \ref{teo2.1}}
Lemma \ref{lem4.5} implies that (5)
of Jakubowski \cite{j96} holds, that is, there exist a dense set $Q$ such that
$(g(S_n(t)), \widetilde{W}_n(t, x))\stackrel{\mathrm
{f.d.d.}}{\longrightarrow
} (g(S(t)), W(t, x)).$
Further, since $g(\cdot)$ is a Lipschitz continuous function and
$S_n(t)$ is uniformly $S$-tight, it follows that
$g(S_n(t))$ is also uniformly $S$-tight. Moreover, for any $x\in
\mathbb{R}$, $\widetilde{W}_n(t, x)$ is a martingale satisfying UT
condition and is $J_1$-tight with limiting law concentrated on $C([0,
1])$, by Remark 4 of Jakubowski \cite{j96}, we see that his condition (6)
is satisfied. Therefore, for $\{g(S_n(t)), \widetilde{W}_n(t, x)\}$,
all the conditions of
Theorem 3 of Jakubowski \cite{j96} are satisfied, as a result of this
theorem, we have
%
\begin{eqnarray}
\label{4.160}\frac{1}{\sqrt{n}}\sum_{i=1}^{n}
g(S_{i-1}/a_n)\zeta_i(x)\stackrel{\mathcal{L}} {
\longrightarrow} \int_0^1 g \bigl(S(t) \bigr)
\, \mathrm{d} W(t, x).
\end{eqnarray}
Thus, (\ref{2.1}) follows from Lemma \ref{lem4.2}. By Lemma \ref{lem4.3}, for (\ref
{2.2}) it suffices to show that
%
\begin{eqnarray}
\label{4.16a} &&U_n(x):=\frac{1}{\sqrt{n}}\sum
_{i=1}^{n}g(S_{i-1}/a_n)
\zeta_i(x)
\nonumber
\\[-8pt]
\\[-8pt]
&&\quad \stackrel{w} {\Longrightarrow}\quad \int_{0}^{1}g
\bigl(S(t) \bigr) \,\mathrm{d}W(t, x)=:U(x),\qquad \mbox{on } D[-A, A].
\nonumber
\end{eqnarray}
The finite-dimension convergence to (\ref{4.16a}) follows from the
Cram\'{e}r--Wold device and (\ref{4.160}). Next, we show
for any $\varepsilon>0$, there exists a $\delta>0$ such that
%
\begin{eqnarray}
\label{4.16b}P \Bigl\{\sup_{|x-y|\leq\delta
}\big|U_n(x)-U_n(y)\big|>
\varepsilon \Bigr\} \rightarrow0.
\end{eqnarray}
This implies that $U_n(x)$ is tight, as a result, we have (\ref{4.16a}).

Since $S_{n}(t)/a_n\stackrel{S}{\Longrightarrow} S(t)$ and $S_n(0)=S(0)=0$,
it follows that
%
\begin{eqnarray}
\label{4.170}\max_{0\leq t\leq1}\big|g \bigl(S_{n}(t)/a_n \bigr)\big|\stackrel {
\mathcal {L}} {\longrightarrow}\max_{0\leq t\leq1} \big|g \bigl(S(t) \bigr)\big|.
\end{eqnarray}

Let $g_\delta(S_i/a_n)=g(S_i/a_n)I(|g(S_i/a_n)|\leq\delta^{-1/4})$ and
$V_n(x)=\frac{1}{\sqrt{n}}\sum_{i=1}^{n}g_\delta(S_{i-1}/a_n)(\zeta_i(x)-\zeta_i(y)).$
Then $V_n(x)$ is a martingale and by Lemma \ref{lem4.1} and condition (A4),
%
\begin{eqnarray}
\label{4.18} &&\mathrm{E} \Bigl[\sup_{y\leq x\leq y+\delta
}\big|V_n(x)\big|
\Bigr]^2
\nonumber
\\
&&\quad \leq {2\over\delta n}\int_{y}^{y+\delta}
\mathrm{E} \Biggl(\sum_{i=1}^{n}
\biggl[g_{\delta} \biggl({S_{i-1}\over a_n} \biggr)
\bigl(\zeta_i(u)- \zeta_i(y) \bigr) \biggr]
\Biggr)^2 \,\mathrm{d}u
\nonumber
\\
&&\qquad {} +{2\delta\over n}\int_{y}^{y+\delta}
\mathrm{E} \Biggl(\sum_{i=1}^{n}
\biggl[g_{\delta} \biggl({S_{i-1}\over a_n} \biggr)
\zeta'_i(u) \biggr] \Biggr)^2\, \mathrm{d}u
\\
&&\quad \leq{2\delta^{-1/2}\over n}\sum_{i=1}^{n}
\int_{y}^{y+\delta
}\mathrm{E} \bigl\{
\zeta_i(u)-\zeta_i(y) \bigr\}^2 \,\mathrm{d}u
+{2 \delta^{-1/2}\delta^2\over n}\sum
_{i=1}^{n}\sup_{y\leq x\leq
y+\delta}\mathrm{E} \bigl\{
\zeta'_i(x) \bigr\}^2
\nonumber
\\
&&\quad \leq{2\delta^{-1/2}\over n}\sum_{i=1}^{n}
\int_{y}^{y+\delta}\int_{y}^{u}
\mathrm{E} \bigl\{\zeta'_i(a) \bigr\}^2 \,
\mathrm{d}a \,\mathrm{d}u+{2\delta^{-1/2}\delta^2\over n}
\sup_{x\in[-A, A]} \mathrm{E} \bigl\{\zeta'_i(x) \bigr
\}^2 \leq C \delta^{3/2}.
\nonumber
\end{eqnarray}
Note that
\begin{eqnarray*}
&&P \Bigl\{\sup_{|x-y|\leq\delta
}\big|U_n(x)-U_n(y)\big|>4
\varepsilon \Bigr\}
\nonumber
\\
&&\quad \leq C\bigl(1+[A/\delta]\bigr)P \Bigl\{\sup_{y\leq x\leq y+\delta
}\big|V_n(x)\big|>
\varepsilon \Bigr\}+P \Bigl\{\max_{1\leq i\leq n}\big|g(S_i/a_n)\big|>
\delta^{-1/4} \Bigr\}.
\nonumber
\end{eqnarray*}
By (\ref{4.170}), (\ref{4.18}) and taking $\delta$ small enough, we
have (\ref{4.16b}) as desired. This completes the proof of Theorem
\ref{teo2.1}.\vspace*{-1pt}
\end{pf*}

\begin{pf*}{Proof of Theorem \ref{teo2.2}}
Note that when $\beta=1$,
$X_i=X_0+\sum_{j=1}^{i}\varepsilon_j$ and $X_0/a_n\stackrel
{p}{\longrightarrow}0$, (\ref{2.4a}) follows directly from (\ref{2.1}).
Next, we show (\ref{2.4b}).

Let $\{u_{ni}\}$ be a constant sequence with $\max_{i}|u_{ni}|=\mathrm{o}(1).$
Along the lines of proof in Lemma~\ref{lem4.2}, we have\vspace*{-2pt}
%
\begin{eqnarray}
\label{4.21} &&\frac{1}{\sqrt{n}}\sum_{i=1}^{n}g
\biggl({S_{i-1}\over a_n} \biggr) \bigl[I(
\varepsilon_i \leq x+u_{ni})-I(\varepsilon_i
\leq x) \bigr]
\nonumber
\\[-1pt]
&&\quad =\frac{1}{\sqrt{n}}\sum_{i=1}^{n}g
\biggl({S_{i-1}\over a_n} \biggr) \bigl[\zeta_i(x+u_{ni})-
\zeta_i(x) \bigr]
\nonumber
\\[-10pt]
\\[-10pt]
&&\qquad {} +\frac{1}{\sqrt{n}}\sum_{i=1}^{n}g
\biggl({S_{i-1}\over a_n} \biggr) \bigl(F(x+u_{ni})-F(x)
\bigr)+\mathrm{o}_p(1)
\nonumber
\\[-1pt]
&&\quad =\frac{1}{\sqrt{n}}\sum_{i=1}^{n}g
\biggl({S_{i-1}\over a_n} \biggr) \bigl(F(x+u_{ni})-F(x)
\bigr)+\mathrm{o}_p(1).
\nonumber\vspace*{-1pt}
\end{eqnarray}

Since
$\max_{1\leq i\leq n}|X_i/a_n|=\mathrm{O}_p(1),$ it follows that when $a_n(\hat
{\beta}-\beta)=\mathrm{o}_p(1)$,
$ \max_{1\leq i\leq n}(\hat{\beta}-\beta)X_i=\mathrm{o}_p(1).$ Thus,
by (\ref{4.21}), we have\vspace*{-1pt}
\begin{eqnarray*}
\frac{1}{\sqrt{n}} \bigl(\hat\alpha_n(x)-\alpha_n(x)
\bigr) 
&=& \frac{1}{\sqrt{n}}\sum
_{i=1}^{n}g(S_{i-1}/a_n) \bigl(F
\bigl(x+(\hat\beta -\beta )X_{i-1} \bigr)-F(x) \bigr)+
\mathrm{o}_p(1).
\nonumber\vspace*{-1pt}
\end{eqnarray*}
This gives (\ref{2.4b}) and completes the proof of Theorem \ref{teo2.2}.\vspace*{-2pt}
\end{pf*}

\begin{pf*}{Proof of Theorem \ref{teo2.3}}
Since $X_i=X_0+\sum_{j=1}^{i}\varepsilon_j$ and $X_0/a_n\stackrel{p}{\longrightarrow
}0$, (\ref{2.5a}) follows from (\ref{2.2}) and the continuous mapping
theorem. Let $u_{ni}$ be given as that in the proof of (\ref{4.21}),
then by Lemma \ref{lem4.3} and a similar argument of (\ref{4.16b}), we have
that under the condition~(A4),\vspace*{-1pt}
\begin{eqnarray*}
\label{4.23} &&\sup_{x\in[-A, A]}\frac{1}{\sqrt{n}}\sum
_{i=1}^{n}g(X_{i-1}/a_n)
\bigl[I( \varepsilon_i\leq x+u_{ni}) \bigr]
\nonumber
\\[-1pt]
&&\quad =\sup_{x\in[-A, A]} \Biggl[\frac{1}{\sqrt{n}}\sum
_{i=1}^{n}g(X_{i-1}/a_n)I(
\varepsilon_i\leq x)+\frac{1}{\sqrt {n}}\sum
_{i=1}^{n}g(X_{i-1}/a_n) \bigl[
\zeta_i(x+u_{ni})-\zeta_i(x) \bigr]
\nonumber
\\[-1pt]
&&\hspace*{60pt} {}+\frac{1}{\sqrt{n}}\sum_{i=1}^{n}g(X_{i-1}/a_n)
\bigl(F(x+u_{ni})-F(x) \bigr) \Biggr]+\mathrm{o}_p(1)
\nonumber
\\[-1pt]
&&\quad =\sup_{x\in[-A, A]} \Biggl[\frac{1}{\sqrt{n}}\sum
_{i=1}^{n}g(X_{i-1}/a_n)I(
\varepsilon_i\leq x)+\frac{1}{\sqrt {n}}\sum
_{i=1}^{n}g(X_{i-1}/a_n)
\bigl(F(x+u_{ni})-F(x) \bigr) \Biggr]
\nonumber
\\[-1pt]
&&\hspace*{21pt} {}+ \mathrm{o}_p(1).
\nonumber
\end{eqnarray*}
As a result, by $\max_{1\leq i\leq n}(\hat{\beta}-1)X_i=\mathrm{o}_p(1)$ and
Taylor's expansion, we have in probability,
%
\begin{eqnarray}
\label{PT2.3a}\sup_{x\in[-A, A]}{\hat\alpha_n(x)\over\sqrt{n}}
&=&\sup_{x\in[-A, A]} \Biggl[{
\alpha_n(x)\over\sqrt{n}}+{1\over \sqrt{n}}f(x) (\hat\beta-1)\sum
_{i=1}^{n}g(X_{i-1}/a_n)X_{i-1}
\Biggr].
\end{eqnarray}
Further, by Theorem 3 of Jakubowski \cite{j96}, it follows that
%
\begin{eqnarray}
\label{PT2.3b}\frac{1}{ n}\sum_{i=1}^{n}
\bigl[g(X_{i-1}/a_n)X_{i-1}/a_n \bigr]
\stackrel{w} {\longrightarrow }\int_0^1 g
\bigl(S(t) \bigr)S(t) \,\mathrm{d}t.
\end{eqnarray}
Combining equations (\ref{2.5a}), (\ref{PT2.3a}), (\ref{PT2.3a})
with Corollary \ref{cor2.1} yields (\ref{2.5b}). Equations (\ref{2.6a}) and
(\ref{2.6b}) follow similarly by noting that when $\hat\beta$ is the
LSE of $\beta$, then
\begin{eqnarray*}
n(\hat\beta-\beta) &=& {1\over2} \Biggl[X_n^2/a_n^2-X_0^2/a_n^2-
\sum_{i=1}^{n}\varepsilon_i^2/a_n^2
\Biggr]\bigg/ \Biggl[\frac{1}{ n}\sum_{i=1}^{n}X_{i-1}^2/a_n^2
\Biggr]
\nonumber
\\
&\stackrel{w} {\longrightarrow}&{1\over2} \bigl(S^2(1)-S^2
\bigr)\Big/\int_0^1 S^2(t) \,
\mathrm{d}t\\
&=:&\int_0^1S(t-) \,\mathrm{d}S(t)\Big/
\int_0^1 S^2(t) \,\mathrm{d}t.
\nonumber
\end{eqnarray*}
The proof of Theorem \ref{teo2.3} is completed.
\end{pf*}

\begin{pf*}{Proof of Theorem \ref{teo3.1}}
Since the proof of the three
cases are similar, we only give the proof under condition (b) in details.
Let $U_{l, i}=\sum_{0\leq j_1<\cdots<j_i}\prod_{s=1}^{i}c_{j_s}\eta_{l-j_s}, U_{l, 0}=1$ and
$L(\widetilde{\varepsilon}_l, x, k)=I(\varepsilon_l\leq x)-\sum_{i=0}^{k}(-i)^{i}F^{(i)}(x)U_{l, i}.$
By Lemma 10 of Wu \cite{w03}, we have that for all $x$,
%
\begin{eqnarray}
\big\Vert \mathcal{P}_1 \bigl(L(\widetilde{\varepsilon}_i,
x, 3) \bigr)\big\Vert &=&\mathrm{O} \Biggl\{|c_{i-1}| \Biggl[|c_{i-1}|+
\Biggl(\sum_{j=i}^{\infty}|c_j|^4
\Biggr)^{1/2}+ \Biggl(\sum_{j=i}^{\infty}|c_j|^2
\Biggr)^{1/2} \Biggr] \Biggr\}
\nonumber
\\[-8pt]
\\[-8pt]
&=&\mathrm{O} \bigl(i^{-2\theta}l^2(i)+i^{-4\theta+3/2}l^3(i)
\bigr).
\nonumber
\end{eqnarray}
Thus, when $\theta>3/4$ or $\theta=3/4$ and $\sum_{i=1}^{\infty
}l^4(i)/i<\infty,$ for all $x$,
\begin{eqnarray*}
\sum_{j=1}^{\infty}\Bigg\Vert \sum
_{i=j}^{\infty}\mathcal {P}_1 \bigl(L(
\widetilde{\varepsilon}_i, x, 3) \bigr)\Bigg\Vert^2&\leq& \sum
_{j=1}^{\infty} \Biggl(\sum
_{i=j}^{\infty}\big\Vert \mathcal {P}_1 \bigl(L(
\widetilde{\varepsilon}_i, x, 3) \bigr)\big\Vert \Biggr)^2\\
& =&
\mathrm{O} \Biggl[\sum_{j=1}^{\infty} \Biggl(\sum
_{i=j}^{\infty}i^{-2\theta}l^2(i)
\Biggr)^2 \Biggr]<\infty.
\nonumber
\end{eqnarray*}
By Theorem 2 of Voln\'{y} \cite{v93}, there exists a martingale difference sequence
$ D_i(x)\in L^2$ and a finite variance
sequence $\{e_i(x)\}$ such that for all $x$,
\[
L(\widetilde{\varepsilon}_i, x, 3)=D_i(x)+e_i(x)-e_{i+1}(x).
\]
Applying (ii) instead of (i) of Lemma \ref{lem4.1} in proving Lemma \ref{lem4.3}, we have that
\[
\sup_{x\in\mathbb{R}}\Bigg|\frac{1}{\sqrt{n}}\sum_{i=1}^{n}g(S_{i-1}/a_n)L(
\widetilde{\varepsilon}_i, x, 3)-{1\over\sqrt{n}}\sum
_{i=1}^{n}g(S_{i-1}/a_n)D_i(x)\Bigg|=
\mathrm{o}_p(1).
\]
Let $g_M(x)=g(x)I(|g(x)|\leq M).$ By (ii) of Lemma \ref{lem4.1}, we have
\begin{eqnarray*}
&&\mathrm{E} \Biggl(\sup_{x\in\mathbb{R}}\frac{1}{\sqrt {n}}\sum
_{i=1}^{n}g_M(S_{i-1}/a_n)D_i(x)
\Biggr)^2
\nonumber
\\
&&\quad \leq{2\over n}\mathrm{E} \int_{\mathbb{R}} \Biggl(\sum
_{i=1}^{n}g_M(S_{i-1}/a_n)D_i(x)
\Biggr)^2 \,\mathrm{d}x+{2\over n}\mathrm{E} \int
_{\mathbb{R}} \Biggl(\sum_{i=1}^{n}g_M(S_{i-1}/a_n)D'_i(x)
\Biggr)^2 \,\mathrm{d}x 
=
\mathrm{O} \bigl(M^2 \bigr).
\nonumber
\end{eqnarray*}
%
As a results, for any positive constants $\varepsilon$ and $\eta$,
there exist a large $M_0$ and a large $N_0$ such that for all $M>M_0$
and $n>N_0$,
%
\begin{eqnarray}
\label{4.29} &&P \Biggl\{\sup_{x\in\mathbb{R}}\frac{1}{ a_n} \Bigg|\sum
_{i=1}^{n}g \biggl({S_{i-1}\over
a_n} \biggr)L(\widetilde {\varepsilon}_i, x, 3) \Bigg|>2
\varepsilon \Biggr\}
\nonumber
\\[-8pt]
\\[-8pt]
&&\quad \leq P \Biggl\{\sup_{x\in\mathbb{R}}\frac{1}{ a_n} \Bigg|\sum
_{i=1}^{n}g_M \biggl({S_{i-1}
\over a_n} \biggr)D_i(x) \Bigg|>\varepsilon \Biggr\}+P \biggl
\{\max_{1\leq i\leq n}\Bigg|g \biggl({S_i\over a_n} \biggr)\Bigg|>M
\biggr\}+\eta \leq3\eta.\qquad \quad
\nonumber
\end{eqnarray}
Note that
%
\begin{eqnarray}
\label{4.31} \frac{1}{ a_n}\sum_{i=1}^{n}g
\biggl({S_{i-1}\over a_n} \biggr)U_{i, 2} &=&
\biggl({{n^{2-2\theta}l^2(n)}\over a_n} \biggr) g
\biggl({S_{n-1}\over a_n} \biggr)\sum
_{i=1}^{n}{U_{i, 2}\over n^{1-2(\theta-1/2)}}
\nonumber
\\[-8pt]
\\[-8pt]
&&{}-{{n^{2-2\theta}l^2(n)}\over a_n}\sum
_{i=1}^{n-1}\sum_{j=1}^{i}
{U_{i, 2}\over{n^{2-2\theta}l^2(n)}} \biggl[g
\biggl({S_{i}\over a_n} \biggr) -g \biggl({S_{i-1}
\over a_n} \biggr) \biggr].\qquad
\nonumber
\end{eqnarray}
From Avram and Taqqu \cite{at87}, it follows that there exists a constant
$C(\theta)$ such that
%
\begin{eqnarray}
\label{4.30}\sum_{i=1}^{[nt]}{U_{i, 2}
\over {n^{2-2\theta
}l^2(n)}}&\stackrel{\mathcal{L}} {
\longrightarrow}&C(\theta)\int_{\mathbb
{R}}\int_{\mathbb{R}}
\int_{0}^{t} \prod_{i=1}^{2}
\bigl[\max(0, v-u_i) \bigr]^{-\theta} \,\mathrm{d}v \,
\mathrm{d}B(u_1) \,\mathrm{d}B(u_2)
\nonumber
\\[-8pt]
\\[-8pt]
&=:&Z_{2, \theta}(t).
\nonumber
\end{eqnarray}
By (\ref{4.30}), the Lipschitz condition of $g(\cdot)$ and an
argument similar to Theorem 3.1 of Ling and Li \cite{ll98}, we have that the
right-hand side of (\ref{4.31}) converges to zero in probability.
Further, by Lemma \ref{lem4.1}, we have
\[
\sup_{x\in\mathbb{R}} \bigl(f^{(k)}(x) \bigr)^2\leq2\int
_{\mathbb
{R}} \bigl(f^{(k)}(x) \bigr)^2 \,
\mathrm{d}x+2\int_{\mathbb{R}} \bigl(f^{(k+1)}(x)
\bigr)^2 \,\mathrm{d}x< \infty\qquad \mbox{for all } k\leq p-1.
\]
Thus, by (\ref{4.29}), we have
$\sup_{x\in\mathbb{R}}\frac{1}{ a_n}\sum_{i=1}^{n}g(S_{i-1}/a_n)[I(\varepsilon_i\leq x)-F(x)+f(x)\varepsilon_i]=\mathrm{o}_p(1).$
Combining this with Lemma \ref{lem4.6} gives that
%
\begin{eqnarray}
\label{4.33}\sup_{x\in\mathbb{R}}\frac{1}{
a_n}\sum
_{i=1}^{n}g(S_{i-1}/a_n)
\bigl[I( \varepsilon_i\leq x)-F(x) \bigr] 
&
\stackrel{ \mathcal{L}} {\longrightarrow}&\sup_{x\in\mathbb
{R}}f(x)\int
_{0}^{1}g \bigl(Z_\theta(t) \bigr) \,
\mathrm{d}Z_\theta(t).
\end{eqnarray}

When $1/2<\theta<3/4$, by (ii) of Theorem 3 in Wu \cite{w03}, we have
%
\begin{eqnarray}
\label{4.25} {1\over n^{2-2\theta}l^2(n)}\sum
_{i=1}^{n} \bigl[I(\varepsilon_i\leq
x)-F(x)+f(x)\varepsilon_i \bigr]\stackrel {w} {\Longrightarrow}f^{\prime}(x)Z_{2, \theta}(1),
\qquad \mbox{on } D(\mathbb{R}).
\end{eqnarray}
Using (\ref{4.25}), we also have (\ref{4.33}). By noting that as $\{
X_t\}$ is a unit root process, $X_t=S_t+X_0.$ This
completes the proof of (\ref{3.1}).

Applying (\ref{4.33}) and arguing as in Theorem \ref{teo2.3}, we have that when
$n(\hat\beta-1)=\mathrm{O}_p(1)$,
\[
\sup_{x\in\mathbb{R}} \Biggl[\frac{1}{ a_n} \bigl(\hat{ \alpha}_n(x)-
\alpha_n(x) \bigr) -\frac{1}{ a_n} \sum
_{i=1}^{n}g(X_{i-1}/a_n) \bigl[F
\bigl(x+(\hat\beta-1)X_{i-1} \bigr)-F(x) \bigr] \Biggr]=
\mathrm{o}_p(1).
\]\eject\noindent
Since $\sup_{x\in\mathbb{R}}f(x)<\infty$ and $\sup_{1\leq i\leq
n}(\hat\beta-1)X_{i}=\mathrm{O}_{p}(a_n/n),$ it follows from Taylor's
expansion and (3.1) that
\begin{eqnarray*}
\sup_{x\in\mathbb{R}}\frac{1}{ a_n} \Biggl[\hat{\alpha}_n(x)
-f(x) (\hat\beta-1)\sum_{i=1}^{n}g(X_{i-1}/a_n)X_{i-1}
\Biggr] 
&\stackrel{\mathcal{L}} {\longrightarrow}&
\sup_{x\in\mathbb
{R}}f(x)\int_{0}^{1}g
\bigl(Z_\theta(t) \bigr) \,\mathrm{d}Z_\theta(t).
\nonumber
\end{eqnarray*}
This gives (\ref{3.2}) and completes the proof of Theorem \ref{teo3.1}.
\end{pf*}

\section*{Acknowledgements}
We would like to thank the Editor, the Associate Editor and two
anonymous referees for helpful comments,
which led to an improved version of this paper. This research was supported
in part by grants from the General Research Fund of HKSAR-RGC-GRF Nos. 400408 and 400410, Collaborative Research Fund of HKSAR-RGC-CRF:
CityU8/CRF/09, Fundamental Research Funds for the Central Universities, ZJNSF (No.~R6090034) and NSFC (Nos. 11171074 and 10801118).
Part of this research was completed when the second author (the
corresponding author) visited CUHK in 2011--2012. Research support from
the Statistics Department
of CUHK is gratefully acknowledged.
%

\printhistory


\begin{thebibliography}{40}

\bibitem{at87}
%
\begin{barticle}[mr]
\bauthor{\bsnm{Avram},~\bfnm{Florin}\binits{F.}} \AND
\bauthor{\bsnm{Taqqu},~\bfnm{Murad~S.}\binits{M.S.}}
(\byear{1987}).
\btitle{Noncentral limit theorems and {A}ppell polynomials}.
\bjournal{Ann. Probab.}
\bvolume{15}
\bpages{767--775}.
\bid{issn={0091-1798}, mr={0885142}}
\bptok{imsref}%
\end{barticle}
%
\endbibitem

\bibitem{at92}
%
\begin{barticle}[mr]
\bauthor{\bsnm{Avram},~\bfnm{Florin}\binits{F.}} \AND
\bauthor{\bsnm{Taqqu},~\bfnm{Murad~S.}\binits{M.S.}}
(\byear{1992}).
\btitle{Weak convergence of sums of moving averages in the {$\alpha$}-stable
domain of attraction}.
\bjournal{Ann. Probab.}
\bvolume{20}
\bpages{483--503}.
\bid{issn={0091-1798}, mr={1143432}}
\bptok{imsref}%
\end{barticle}
%
\endbibitem

\bibitem{b03}
%
\begin{barticle}[auto:STB|2012/08/01|11:33:29]
\bauthor{\bsnm{Bai},~\bfnm{J.~S.}\binits{J.S.}}
(\byear{2003}).
\btitle{Testing parametric conditional distributions of dynamic models}.
\bjournal{The Review of Economics and Statistics}
\bvolume{85}
\bpages{531--549}.
\bptok{imsref}%
\end{barticle}
%
\endbibitem

\bibitem{b96}
%
\begin{barticle}[mr]
\bauthor{\bsnm{Baillie},~\bfnm{Richard~T.}\binits{R.T.}}
(\byear{1996}).
\btitle{Long memory processes and fractional integration in econometrics}.
\bjournal{J. Econometrics}
\bvolume{73}
\bpages{5--59}.
\bid{doi={10.1016/0304-4076(95)01732-1}, issn={0304-4076}, mr={1410000}}
\bptok{imsref}%
\end{barticle}
%
\endbibitem

\bibitem{bw71}
%
\begin{barticle}[mr]
\bauthor{\bsnm{Bickel},~\bfnm{P.~J.}\binits{P.J.}} \AND
\bauthor{\bsnm{Wichura},~\bfnm{M.~J.}\binits{M.J.}}
(\byear{1971}).
\btitle{Convergence criteria for multiparameter stochastic processes
and some
applications}.
\bjournal{Ann. Math. Statist.}
\bvolume{42}
\bpages{1656--1670}.
\bid{issn={0003-4851}, mr={0383482}}
\bptok{imsref}%
\end{barticle}
%
\endbibitem

\bibitem{b68}
%
\begin{bbook}[mr]
\bauthor{\bsnm{Billingsley},~\bfnm{Patrick}\binits{P.}}
(\byear{1968}).
\btitle{Convergence of Probability Measures}.
\baddress{New York}: \bpublisher{Wiley}.
\bid{mr={0233396}}
\bptok{imsref}%
\end{bbook}
%
\endbibitem

\bibitem{cl08}
%
\begin{barticle}[mr]
\bauthor{\bsnm{Chan},~\bfnm{Ngai~Hang}\binits{N.H.}} \AND
\bauthor{\bsnm{Ling},~\bfnm{Shiqing}\binits{S.}}
(\byear{2008}).
\btitle{Residual empirical processes for long and short memory time series}.
\bjournal{Ann. Statist.}
\bvolume{36}
\bpages{2453--2470}.
\bid{doi={10.1214/07-AOS543}, issn={0090-5364}, mr={2458194}}
\bptok{imsref}%
\end{barticle}
%
\endbibitem

\bibitem{cz09b}
%
\begin{barticle}[mr]
\bauthor{\bsnm{Chan},~\bfnm{Ngai~Hang}\binits{N.H.}} \AND
\bauthor{\bsnm{Zhang},~\bfnm{Rong-Mao}\binits{R.M.}}
(\byear{2009}).
\btitle{Inference for nearly nonstationary processes under strong dependence
with infinite variance}.
\bjournal{Statist. Sinica}
\bvolume{19}
\bpages{925--947}.
\bid{issn={1017-0405}, mr={2536137}}
\bptok{imsref}%
\end{barticle}
%
\endbibitem

\bibitem{cz09a}
%
\begin{barticle}[mr]
\bauthor{\bsnm{Chan},~\bfnm{Ngai~Hang}\binits{N.H.}} \AND
\bauthor{\bsnm{Zhang},~\bfnm{Rong-Mao}\binits{R.M.}}
(\byear{2009}).
\btitle{Quantile inference for near-integrated autoregressive time
series under
infinite variance and strong dependence}.
\bjournal{Stochastic Process. Appl.}
\bvolume{119}
\bpages{4124--4148}.
\bid{doi={10.1016/j.spa.2009.09.010}, issn={0304-4149}, mr={2565561}}
\bptok{imsref}%
\end{barticle}
%
\endbibitem

\bibitem{cz10}
%
\begin{barticle}[mr]
\bauthor{\bsnm{Chan},~\bfnm{Ngai~Hang}\binits{N.H.}} \AND
\bauthor{\bsnm{Zhang},~\bfnm{Rong-Mao}\binits{R.M.}}
(\byear{2010}).
\btitle{Inference for unit-root models with infinite variance {GARCH} errors}.
\bjournal{Statist. Sinica}
\bvolume{20}
\bpages{1363--1393}.
\bid{issn={1017-0405}, mr={2777329}}
\bptok{imsref}%
\end{barticle}
%
\endbibitem

\bibitem{dr86}
%
\begin{barticle}[mr]
\bauthor{\bsnm{Davis},~\bfnm{Richard}\binits{R.}} \AND
\bauthor{\bsnm{Resnick},~\bfnm{Sidney}\binits{S.}}
(\byear{1986}).
\btitle{Limit theory for the sample covariance and correlation
functions of
moving averages}.
\bjournal{Ann. Statist.}
\bvolume{14}
\bpages{533--558}.
\bid{doi={10.1214/aos/1176349937}, issn={0090-5364}, mr={0840513}}
\bptok{imsref}%
\end{barticle}
%
\endbibitem

\bibitem{dm03}
%
\begin{barticle}[mr]
\bauthor{\bsnm{Dedecker},~\bfnm{J{\'e}r{\^o}me}\binits{J.}} \AND
\bauthor{\bsnm{Merlev{\`e}de},~\bfnm{Florence}\binits{F.}}
(\byear{2003}).
\btitle{The conditional central limit theorem in {H}ilbert spaces}.
\bjournal{Stochastic Process. Appl.}
\bvolume{108}
\bpages{229--262}.
\bid{doi={10.1016/j.spa.2003.07.004}, issn={0304-4149}, mr={2019054}}
\bptok{imsref}%
\end{barticle}
%
\endbibitem

\bibitem{bdg07}
%
\begin{bbook}[mr]
\bauthor{\bparticle{del} \bsnm{Barrio},~\bfnm{Eustasio}\binits{E.}},
\bauthor{\bsnm{Deheuvels},~\bfnm{Paul}\binits{P.}} \AND
\bauthor{\bparticle{van~de} \bsnm{Geer},~\bfnm{Sara}\binits{S.}}
(\byear{2007}).
\btitle{Lectures on Empirical Processes: Theory and Statistical Applications}.
\bseries{EMS Series of Lectures in Mathematics}.
\baddress{Z\"urich}: \bpublisher{Eur. Math. Soc.}
\bnote{With a preface by Juan A. Cuesta
Albertos and Carlos Matr{\'a}n}.
\bid{doi={10.4171/027}, mr={2284824}}
\bptok{imsref}%
\end{bbook}
%
\endbibitem

\bibitem{e06}
%
\begin{barticle}[mr]
\bauthor{\bsnm{Escanciano},~\bfnm{J.~Carlos}\binits{J.C.}}
(\byear{2006}).
\btitle{Goodness-of-fit tests for linear and nonlinear time series models}.
\bjournal{J. Amer. Statist. Assoc.}
\bvolume{101}
\bpages{531--541}.
\bid{doi={10.1198/016214505000001050}, issn={0162-1459}, mr={2256173}}
\bptok{imsref}%
\end{barticle}
%
\endbibitem

\bibitem{e07}
%
\begin{barticle}[mr]
\bauthor{\bsnm{Escanciano},~\bfnm{J.~Carlos}\binits{J.C.}}
(\byear{2007}).
\btitle{Weak convergence of non-stationary multivariate marked
processes with
applications to martingale testing}.
\bjournal{J. Multivariate Anal.}
\bvolume{98}
\bpages{1321--1336}.
\bid{doi={10.1016/j.jmva.2007.03.004}, issn={0047-259X}, mr={2364121}}
\bptok{imsref}%
\end{barticle}
%
\endbibitem

\bibitem{e10}
%
\begin{barticle}[mr]
\bauthor{\bsnm{Escanciano},~\bfnm{J.~Carlos}\binits{J.C.}}
(\byear{2010}).
\btitle{Asymptotic distribution-free diagnostic tests for
heteroskedastic time
series models}.
\bjournal{Econometric Theory}
\bvolume{26}
\bpages{744--773}.
\bid{doi={10.1017/S0266466609990090}, issn={0266-4666}, mr={2646478}}
\bptok{imsref}%
\end{barticle}
%
\endbibitem


\bibitem{gr76}
%
\begin{bbook}[mr]
\bauthor{\bsnm{Gaenssler},~\bfnm{Peter}\binits{P.}} \AND
\bauthor{\bsnm{R{\'e}v{\'e}sz},~\bfnm{P{\'a}l}\binits{P.}}
(\byear{1976}).
\btitle{Empirical Distributions and Processes}.
\bseries{Lecture Notes in Mathematics} \bvolume{566}.
\baddress{Berlin}: \bpublisher{Springer}.
\bid{mr={0428364}}
\bptok{imsref}%
\end{bbook}
%
\endbibitem



\bibitem{gk54}
%
\begin{bbook}[mr]
\bauthor{\bsnm{Gnedenko},~\bfnm{B.~V.}\binits{B.V.}} \AND
\bauthor{\bsnm{Kolmogorov},~\bfnm{A.~N.}\binits{A.N.}}
(\byear{1954}).
\btitle{Limit Distributions for Sums of Independent Random Variables}.
\baddress{Cambridge, MA}: \bpublisher{Addison-Wesley}.
\bnote{Translated and annotated by K.L. Chung. With an Appendix by J.L.
Doob}.
\bid{mr={0062975}}
\bptok{imsref}%
\end{bbook}
%
\endbibitem

\bibitem{hh80}
%
\begin{bbook}[mr]
\bauthor{\bsnm{Hall},~\bfnm{P.}\binits{P.}} \AND
\bauthor{\bsnm{Heyde},~\bfnm{C.~C.}\binits{C.C.}}
(\byear{1980}).
\btitle{Martingale Limit Theory and Its Application (Probability and Mathematical Statistics)}.
\baddress{New York}: \bpublisher{Academic Press [Harcourt Brace Jovanovich
Publishers]}.
\bid{mr={0624435}}
\bptok{imsref}%
\end{bbook}
%
\endbibitem

\bibitem{hh96}
%
\begin{barticle}[mr]
\bauthor{\bsnm{Ho},~\bfnm{Hwai-Chung}\binits{H.C.}} \AND
\bauthor{\bsnm{Hsing},~\bfnm{Tailen}\binits{T.}}
(\byear{1996}).
\btitle{On the asymptotic expansion of the empirical process of long-memory
moving averages}.
\bjournal{Ann. Statist.}
\bvolume{24}
\bpages{992--1024}.
\bid{doi={10.1214/aos/1032526953}, issn={0090-5364}, mr={1401834}}
\bptok{imsref}%
\end{barticle}
%
\endbibitem

\bibitem{hl03}
%
\begin{barticle}[mr]
\bauthor{\bsnm{Hong},~\bfnm{Yongmiao}\binits{Y.}} \AND
\bauthor{\bsnm{Lee},~\bfnm{Tae-Hwy}\binits{T.H.}}
(\byear{2003}).
\btitle{Diagnostic checking for the adequacy of nonlinear time series models}.
\bjournal{Econometric Theory}
\bvolume{19}
\bpages{1065--1121}.
\bid{doi={10.1017/S0266466603196089}, issn={0266-4666}, mr={2015977}}
\bptok{imsref}%
\end{barticle}
%
\endbibitem

\bibitem{j96}
%
\begin{barticle}[mr]
\bauthor{\bsnm{Jakubowski},~\bfnm{Adam}\binits{A.}}
(\byear{1996}).
\btitle{Convergence in various topologies for stochastic integrals
driven by
semimartingales}.
\bjournal{Ann. Probab.}
\bvolume{24}
\bpages{2141--2153}.
\bid{doi={10.1214/aop/1041903222}, issn={0091-1798}, mr={1415245}}
\bptok{imsref}%
\end{barticle}
%
\endbibitem

\bibitem{j97}
%
\begin{barticle}[mr]
\bauthor{\bsnm{Jakubowski},~\bfnm{Adam}\binits{A.}}
(\byear{1997}).
\btitle{A non-{S}korohod topology on the {S}korohod space}.
\bjournal{Electron. J. Probab.}
\bvolume{2}
\bpages{21 (electronic)}.
\bid{doi={10.1214/EJP.v2-18}, issn={1083-6489}, mr={1475862}}
\bptok{imsref}%
\end{barticle}
%
\endbibitem

\bibitem{k73}
%
\begin{barticle}[mr]
\bauthor{\bsnm{Kesten},~\bfnm{Harry}\binits{H.}}
(\byear{1973}).
\btitle{Random difference equations and renewal theory for products of random
matrices}.
\bjournal{Acta Math.}
\bvolume{131}
\bpages{207--248}.
\bid{issn={0001-5962}, mr={0440724}}
\bptok{imsref}%
\end{barticle}
%
\endbibitem

\bibitem{k91}
%
\begin{barticle}[mr]
\bauthor{\bsnm{Knight},~\bfnm{Keith}\binits{K.}}
(\byear{1991}).
\btitle{Limit theory for {$M$}-estimates in an integrated infinite variance
process}.
\bjournal{Econometric Theory}
\bvolume{7}
\bpages{200--212}.
\bid{doi={10.1017/S0266466600004400}, issn={0266-4666}, mr={1128412}}
\bptok{imsref}%
\end{barticle}
%
\endbibitem

\bibitem{kl06}
%
\begin{barticle}[mr]
\bauthor{\bsnm{Koul},~\bfnm{Hira~L.}\binits{H.L.}} \AND
\bauthor{\bsnm{Ling},~\bfnm{Shiqing}\binits{S.}}
(\byear{2006}).
\btitle{Fitting an error distribution in some heteroscedastic time series
models}.
\bjournal{Ann. Statist.}
\bvolume{34}
\bpages{994--1012}.
\bid{doi={10.1214/009053606000000191}, issn={0090-5364}, mr={2283401}}
\bptok{imsref}%
\end{barticle}
%
\endbibitem

\bibitem{kp91}
%
\begin{barticle}[mr]
\bauthor{\bsnm{Kurtz},~\bfnm{Thomas~G.}\binits{T.G.}} \AND
\bauthor{\bsnm{Protter},~\bfnm{Philip}\binits{P.}}
(\byear{1991}).
\btitle{Weak limit theorems for stochastic integrals and stochastic
differential equations}.
\bjournal{Ann. Probab.}
\bvolume{19}
\bpages{1035--1070}.
\bid{issn={0091-1798}, mr={1112406}}
\bptok{imsref}%
\end{barticle}
%
\endbibitem

\bibitem{ll98}
%
\begin{barticle}[mr]
\bauthor{\bsnm{Ling},~\bfnm{Shiqing}\binits{S.}} \AND
\bauthor{\bsnm{Li},~\bfnm{W.~K.}\binits{W.K.}}
(\byear{1998}).
\btitle{Limiting distributions of maximum likelihood estimators for unstable
autoregressive moving-average time series with general autoregressive
heteroscedastic errors}.
\bjournal{Ann. Statist.}
\bvolume{26}
\bpages{84--125}.
\bid{doi={10.1214/aos/1030563979}, issn={0090-5364}, mr={1611800}}
\bptok{imsref}%
\end{barticle}
%
\endbibitem

\bibitem{m68}
%
\begin{barticle}[mr]
\bauthor{\bsnm{Marcus},~\bfnm{M.~B.}\binits{M.B.}}
(\byear{1968}).
\btitle{H\"older conditions for {G}aussian processes with stationary
increments}.
\bjournal{Trans. Amer. Math. Soc.}
\bvolume{134}
\bpages{29--52}.
\bid{issn={0002-9947}, mr={0230368}}
\bptok{imsref}%
\end{barticle}
%
\endbibitem

\bibitem{mn00}
%
\begin{barticle}[mr]
\bauthor{\bsnm{Mikosch},~\bfnm{Thomas}\binits{T.}} \AND
\bauthor{\bsnm{Norvai{\v{s}}a},~\bfnm{Rimas}\binits{R.}}
(\byear{2000}).
\btitle{Stochastic integral equations without probability}.
\bjournal{Bernoulli}
\bvolume{6}
\bpages{401--434}.
\bid{doi={10.2307/3318668}, issn={1350-7265}, mr={1762553}}
\bptok{imsref}%
\end{barticle}
%
\endbibitem

\bibitem{rg79}
%
\begin{barticle}[mr]
\bauthor{\bsnm{Resnick},~\bfnm{Sidney}\binits{S.}} \AND
\bauthor{\bsnm{Greenwood},~\bfnm{Priscilla}\binits{P.}}
(\byear{1979}).
\btitle{A bivariate stable characterization and domains of attraction}.
\bjournal{J. Multivariate Anal.}
\bvolume{9}
\bpages{206--221}.
\bid{doi={10.1016/0047-259X(79)90079-4}, issn={0047-259X}, mr={0538402}}
\bptok{imsref}%
\end{barticle}
%
\endbibitem

\bibitem{sxz08}
%
\begin{barticle}[mr]
\bauthor{\bsnm{Stute},~\bfnm{W.}\binits{W.}},
\bauthor{\bsnm{Xu},~\bfnm{W.~L.}\binits{W.L.}} \AND
\bauthor{\bsnm{Zhu},~\bfnm{L.~X.}\binits{L.X.}}
(\byear{2008}).
\btitle{Model diagnosis for parametric regression in high-dimensional spaces}.
\bjournal{Biometrika}
\bvolume{95}
\bpages{451--467}.
\bid{doi={10.1093/biomet/asm095}, issn={0006-3444}, mr={2521592}}
\bptok{imsref}%
\end{barticle}
%
\endbibitem


\bibitem{tk07}
%
\begin{bbook}[mr]
\bauthor{\bsnm{Teyssi{\`e}re},~\bfnm{Gilles}\binits{G.}} \AND
\bauthor{\bsnm{Kirman},~\bfnm{Alan P.}\binits{A. P.}}
(\byear{2007}).
\btitle{Long Memory in Economics}.
\baddress{Berlin}: \bpublisher{Springer}.
\bid{doi={10.1007/978-3-540-34625-8}, mr={2263582}}
\bptok{imsref}%
\end{bbook}
%
\endbibitem


\bibitem{v93}
%
\begin{barticle}[mr]
\bauthor{\bsnm{Voln{\'y}},~\bfnm{Dalibor}\binits{D.}}
(\byear{1993}).
\btitle{Approximating martingales and the central limit theorem for strictly
stationary processes}.
\bjournal{Stochastic Process. Appl.}
\bvolume{44}
\bpages{41--74}.
\bid{doi={10.1016/0304-4149(93)90037-5}, issn={0304-4149}, mr={1198662}}
\bptok{imsref}%
\end{barticle}
%
\endbibitem

\bibitem{w90}
%
\begin{barticle}[mr]
\bauthor{\bsnm{Wooldridge},~\bfnm{Jeffrey~M.}\binits{J.M.}}
(\byear{1990}).
\btitle{A unified approach to robust, regression-based specification tests}.
\bjournal{Econometric Theory}
\bvolume{6}
\bpages{17--43}.
\bid{doi={10.1017/S0266466600004898}, issn={0266-4666}, mr={1059144}}
\bptok{imsref}%
\end{barticle}
%
\endbibitem

\bibitem{w03}
%
\begin{barticle}[mr]
\bauthor{\bsnm{Wu},~\bfnm{Wei~Biao}\binits{W.B.}}
(\byear{2003}).
\btitle{Empirical processes of long-memory sequences}.
\bjournal{Bernoulli}
\bvolume{9}
\bpages{809--831}.
\bid{doi={10.3150/bj/1066418879}, issn={1350-7265}, mr={2047687}}
\bptok{imsref}%
\end{barticle}
%
\endbibitem

\bibitem{w07}
%
\begin{barticle}[mr]
\bauthor{\bsnm{Wu},~\bfnm{Wei~Biao}\binits{W.B.}}
(\byear{2007}).
\btitle{Strong invariance principles for dependent random variables}.
\bjournal{Ann. Probab.}
\bvolume{35}
\bpages{2294--2320}.
\bid{doi={10.1214/009117907000000060}, issn={0091-1798}, mr={2353389}}
\bptok{imsref}%
\end{barticle}
%
\endbibitem

\bibitem{y36}
%
\begin{barticle}[mr]
\bauthor{\bsnm{Young},~\bfnm{L.~C.}\binits{L.C.}}
(\byear{1936}).
\btitle{An inequality of the {H}\"older type, connected with {S}tieltjes
integration}.
\bjournal{Acta Math.}
\bvolume{67}
\bpages{251--282}.
\bid{doi={10.1007/BF02401743}, issn={0001-5962}, mr={1555421}}
\bptok{imsref}%
\end{barticle}
%
\endbibitem

\bibitem{zy08}
%
\begin{barticle}[mr]
\bauthor{\bsnm{Zou},~\bfnm{Hui}\binits{H.}} \AND
\bauthor{\bsnm{Yuan},~\bfnm{Ming}\binits{M.}}
(\byear{2008}).
\btitle{Composite quantile regression and the oracle model selection theory}.
\bjournal{Ann. Statist.}
\bvolume{36}
\bpages{1108--1126}.
\bid{doi={10.1214/07-AOS507}, issn={0090-5364}, mr={2418651}}
\bptok{imsref}%
\end{barticle}
%
\endbibitem

\end{thebibliography}
\end{document}